\documentclass[moor]{informs1}              

\usepackage{natbib}
 \NatBibNumeric
 \def\bibfont{\small}%
 \bibpunct[, ]{[}{]}{,}{n}{}{,}%

\usepackage[colorlinks=true,breaklinks=true,bookmarks=true,urlcolor=blue,
     citecolor=blue,linkcolor=blue,bookmarksopen=false,draft=false]{hyperref}

\def\EMAIL#1{\href{mailto:#1}{#1}}
\def\URL#1{\href{#1}{#1}}         

\TheoremsNumberedThrough     

\EquationsNumberedThrough    

\usepackage{convex}
\usepackage{algorithm}
\usepackage{algpseudocode}
\usepackage{enumitem}

\begin{document}


 \RUNAUTHOR{Hauenstein, Mohammad-Nezhad, Tang, and Terlaky} 

\RUNTITLE{On computing the nonlinearity interval in parametric SDO}

 \TITLE{On computing the nonlinearity interval in parametric semidefinite optimization}

\ARTICLEAUTHORS{%
\AUTHOR{Jonathan D. Hauenstein}
\AFF{Department of Applied and Computational Mathematics and Statistics, University of Notre Dame, \EMAIL{hauenstein@nd.edu} \URL{}}
\AUTHOR{Ali Mohammad-Nezhad}
\AFF{Department of Mathematics, Purdue University, \EMAIL{mohamm42@purdue.edu} \URL{}}
\AUTHOR{Tingting Tang}
\AFF{Department of Applied and Computational Mathematics and Statistics, University of Notre Dame, \EMAIL{ttang@nd.edu} \URL{}}
\AUTHOR{Tam\'as Terlaky}
\AFF{Department of Industrial and Systems Engineering, Lehigh University, \EMAIL{terlaky@lehigh.edu} \URL{}}
} 

\ABSTRACT{%
This paper revisits the parametric analysis of semidefinite optimization problems with respect to the perturbation of the objective function along a fixed direction. We review the notions of invariancy set, nonlinearity interval, and transition point of the optimal partition, and we investigate their characterizations. We show that the set of transition points is finite and the continuity of the optimal set mapping, on the basis of Painlev\'e-Kuratowski set convergence, might fail on a nonlinearity interval. Under a local nonsingularity condition, we then develop a methodology, stemming from numerical algebraic geometry, to efficiently compute nonlinearity intervals and transition points of the optimal partition. Finally, we support the theoretical results by applying our procedure to some numerical examples.
}%


\KEYWORDS{ Parametric semidefinite optimization; Optimal partition; Nonlinearity interval; Numerical algebraic geometry}
\MSCCLASS{Primary: 90C22; Secondary: 90C31, 90C51}

\maketitle

\section{Introduction}\label{intro}
Let $\mathbb{S}^n$ be the vector space of $n \times n$ symmetric matrices. Consider a parametric semidefinite optimization (SDO) problem
 \begin{align*}
&(\mathrm{P_{\epsilon}}) \qquad \inf_{X \in \mathbb{S}^n} \big \{\langle C + \epsilon \bar{C},  X \rangle : \langle A^{i} , X \rangle=b_i, \quad  i=1,\ldots, m, \ X \succeq 0 \big \},\\
&(\mathrm{D_{\epsilon}}) \qquad \sup_{(y,S) \in \mathbb{R}^m \times \mathbb{S}^n} \ \bigg \{b^Ty  : \sum_{i=1}^m y_i A^i+S=C + \epsilon \bar{C}, \ S \succeq 0 \bigg \},
 \end{align*}
where $C,A^i \in \mathbb{S}^n$ for $i=1,\ldots, m$, $b \in \mathbb{R}^m$, $\bar{C} \in \mathbb{S}^n$ is a fixed direction, the inner product is defined as $\langle C,X \rangle\!:=\trace(CX)$, and $X \succeq 0$ means that the matrix $X$ is symmetric and positive semidefinite. Let $v(\epsilon)\in \mathbb{R} \cup \{-\infty,\infty\}$ 
denote the optimal value of $(\mathrm{P_{\epsilon}})$. This yields a function $v:\mathbb{R} \rightarrow \mathbb{R} \cup \{-\infty,\infty\}$ which is the so-called \textit{optimal value function}. Let $\mathcal{E}\!:=\{\epsilon \in \mathbb{R}:v(\epsilon) > -\infty\}$ be the domain of~$v(\epsilon)$.

\vspace{5px}
\noindent
The primal and dual optimal set mappings on $\mathcal{E}$ are defined as
\begin{align*}
&\mathcal{P}^*:\epsilon \mapsto \big\{X : \langle C + \epsilon \bar{C}, X \rangle = v(\epsilon), \ X \in \mathcal{P}(\epsilon)\big\}, \\
&\mathcal{D}^*:\epsilon \mapsto \big \{(y,S) : b^T y = v(\epsilon), \ (y,S) \in \mathcal{D}(\epsilon)\big \},
\end{align*}
where $\mathcal{P}$ and $\mathcal{D}$ denote the primal and dual feasible set mappings:
\begin{align*}
&\mathcal{P}: \epsilon \mapsto \big\{X : \langle A^{i} , X \rangle=b_i, \  i=1,\ldots, m, \ X \succeq 0\big\}, \\
&\mathcal{D}: \epsilon \mapsto\!\bigg\{(y,S) : \sum_{i=1}^m y_i A^i+S=C + \epsilon \bar{C}, \ S \succeq 0 \bigg\}.
\end{align*}
Note that $\mathcal{P}^*(\epsilon)$ or $\mathcal{D}^*(\epsilon)$ might be empty for some $\epsilon \in \mathcal{E}$. To avoid trivialities, we make the following assumptions throughout this paper: 
\begin{assumption}\label{IND}
The coefficient matrices $A^i$ for $i=1,\ldots, m$ are linearly independent.
\end{assumption}
\begin{assumption}\label{IPC}
The interior point condition holds for both $(\mathrm{P_{\epsilon}})$ and $(\mathrm{D_{\epsilon}})$ at $\epsilon = 0$, i.e., there exists a feasible $\big(X^{\circ}(0),y^{\circ}(0),S^{\circ}(0)\big) \in \mathcal{P}(0) \times \mathcal{D}(0)$ such that $X^{\circ}(0),S^{\circ}(0) \succ 0$, where $\succ 0$ means positive definite.
 \end{assumption}
We may assume Assumption~\ref{IPC} without loss of generality. In fact, the interior point condition is standard in the literature of conic optimization, and it always holds for a self-dual homogeneous embedding form of an SDO problem~\cite{Kl97,KRT98}. Assumption~\ref{IPC} implies that $\mathcal{E}$ is nonempty and non-singleton~\cite[Theorem~4.1]{T01}, and that $v(\epsilon)$ is proper and concave on $\mathcal{E}$. The proof is analogous to~\cite[Theorem~11]{BRT97}, where the objective function is linear. The concavity of $v(\epsilon)$ yields that $\mathcal{E}$ is a closed, possibly unbounded, interval, see e.g.,~\cite[Theorem~8]{BRT97} and that $v(\epsilon)$ is continuous on $\interior(\mathcal{E})$~\cite[Corollary 2.109]{BS00}, where $\interior(\cdot)$ denotes the interior of a set.

\begin{remark}
By~\cite[Lemma~3.1]{GS99} and a theorem of the alternative~\cite[Lemma~12.6]{CSW13}, Assumptions~\ref{IND} and~\ref{IPC} imply that a strictly feasible solution $\big(X^{\circ}(\epsilon),y^{\circ}(\epsilon),S^{\circ}(\epsilon)\big)$ exists at every $\epsilon \in \interior(\mathcal{E})$. \halmos
\end{remark}

\vspace{5px}
\noindent
Hence, for all $\epsilon\in\interior({\mathcal{E}})$,
Assumptions~\ref{IND} and~\ref{IPC} ensure that strong duality holds and that the optimal sets $\mathcal{P}^*(\epsilon)$ and $\mathcal{D}^*(\epsilon)$ are nonempty and compact~\cite[Corollary~4.2]{T01}. In this paper, by strong duality we mean that the optimal values of $(\mathrm{P_{\epsilon}})$ and $(\mathrm{D_{\epsilon}})$ are both attained and the duality gap is zero. In particular, the optimality conditions for $(\mathrm{P_{\epsilon}})$ and $(\mathrm{D_{\epsilon}})$ can be written as 
\begin{equation}\label{KKT_param_SDO}
\begin{aligned}
\langle A^i, X \rangle &= b_i, & i&=1,\ldots, m,\\
\sum_{i=1}^m y_iA^i + S &= C + \epsilon \bar{C},\\
XS&= 0,\\
X,S &\succeq 0,
\end{aligned}
\end{equation} 
where $XS= 0$ denotes the complementarity condition. Furthermore, Assumption~\ref{IPC} guarantees the existence of a 
so-called maximally complementary optimal solution for every $\epsilon \in \interior(\mathcal{E})$.

\begin{definition}
For any fixed $\epsilon\in\interior(\mathcal{E})$, an optimal solution $\big(X^*(\epsilon),y^*(\epsilon),S^*(\epsilon)\big)$ is called \textit{maximally complementary} if 
\begin{align*}
X^*(\epsilon) \in \ri\!\big(\mathcal{P}^*(\epsilon)\big) \quad \text{and} \quad \big(y^*(\epsilon),S^*(\epsilon)\big) \in \ri\!\big(\mathcal{D}^*(\epsilon)\big),
\end{align*}
where $\ri(\cdot)$ denotes the relative interior of a set. A maximally complementary optimal solution $\big(X^*(\epsilon),y^*(\epsilon),S^*(\epsilon)\big)$ is called \textit{strictly complementary} if $X^*(\epsilon)+S^*(\epsilon) \succ 0$.
\end{definition}
\vspace{5px}
\noindent
For a given $\epsilon \in \interior(\mathcal{E})$, unless stated otherwise, $\big(X^*(\epsilon),y^*(\epsilon),S^*(\epsilon)\big)$ denotes a maximally complementary optimal solution. Notice that $\rank\big(X^*(\epsilon)\big)+\rank\big(S^*(\epsilon)\big)$ is maximal on $\mathcal{P}^*(\epsilon) \times \mathcal{D}^*(\epsilon)$, see e.g.,~\cite[Lemma 2.3]{Kl02}. Even though a strictly complementary optimal solution may fail to exist, a maximally complementary optimal solution always exists under Assumption~\ref{IPC}.

\vspace{5px}
\noindent
In practice, given a fixed $\epsilon$, $(\mathrm{P_{\epsilon}})$ and $(\mathrm{D_{\epsilon}})$ can be efficiently solved using a primal-dual path-following interior point method (IPM), see~\cite{NN94}. A primal-dual path following IPM generates a sequence of solutions whose accumulation points are maximally complementary optimal solutions~\cite{H02}.

\subsection{Optimal partition}
For SDO, the optimal partition information can be leveraged to establish sensitivity analysis results. The optimal partition provides a characterization of the optimal set, and it is uniquely defined for any instance of an SDO problem which satisfies strong duality~\cite{Kl02}. For a fixed $\epsilon \in \interior(\mathcal{E})$, let $\big(X^*(\epsilon),y^*(\epsilon),S^*(\epsilon)\big) \in \ri\big(\mathcal{P}^*(\epsilon) \times \mathcal{D}^*(\epsilon)\big)$ be a maximally complementary optimal solution, and let $\mathcal{B}(\epsilon)\!:=\mathcal{R}\big(X^*(\epsilon)\big)$, $\mathcal{N}(\epsilon)\!:=\mathcal{R}\big(S^*(\epsilon)\big)$, and $\mathcal{T}(\epsilon)\!:=\big(\mathcal{R}\big(X^*(\epsilon)\big) + \mathcal{R}\big(S^*(\epsilon)\big)\big)^{\perp}$, where~$\mathcal{R}(\cdot)$ is the column space and $\perp$ denotes the orthogonal complement of a subspace. Then the 3-tuple $\big(\mathcal{B}(\epsilon), \mathcal{T}(\epsilon),\mathcal{N}(\epsilon)\big)$ is called the \textit{optimal partition} of $(\mathrm{P_{\epsilon}})$ and $(\mathrm{D_{\epsilon}})$. Note that the subspaces $\mathcal{R}\big(X^*(\epsilon)\big)$ and $\mathcal{R}\big(S^*(\epsilon)\big)$ are orthogonal by the complementarity condition in~\eqref{KKT_param_SDO}.
Further, the optimal partition $\big(\mathcal{B}(\epsilon), \mathcal{T}(\epsilon),\mathcal{N}(\epsilon)\big)$
is independent of the choice of a maximally complementary optimal solution~\cite[Lemma 2.3(i)]{Kl02}.

\subsection{Related work}
Sensitivity analysis along a fixed direction has been extensively studied
in optimization theory and was originally introduced for linear optimization (LO) and linearly constrained quadratic optimization (LCQO) problems in~\cite{AM92,BJRT96,JRT93}. Sensitivity analysis of nonlinear optimization problems was studied by Fiacco~\cite{Fi76} and Fiacco and McCormick~\cite{FM90} using the implicit function theorem~\cite[Theorem~10.2.1]{D60}. Their analysis was based on linear independence constraint qualification, second-order sufficient condition, and the strict complementarity condition. Furthermore, Fiacco~\cite{Fi76} showed how to compute/approximate the partial derivatives of a locally optimal solution. Robinson~\cite{R82} removed the reliance on the strict complementarity condition by imposing a strong second-order sufficient condition. Kojima~\cite{K80} removed the dependence on the strict complementarity condition by invoking the degree theory of a continuous map, see e.g.,~\cite{OR70}. A comprehensive treatment of directional and differential stability of nonlinear conic optimization problems is given by Bonnans and Shapiro~\cite{BS98,BS00}, see also~\cite{BR2005,S97}. The reader is referred to~\cite{Fi83} for a survey of classical results.

\vspace{5px}
\noindent
 The study of sensitivity analysis based on the optimal partition approach was initiated by Adler and Monteiro~\cite{AM92} and Jansen et al.~\cite{JRT93} for LO and then extended to LCQO, SDO, and linear conic optimization by Berkelaar et al.~\cite{BJRT96}, Goldfarb and Scheinberg~\cite{GS99}, and Yildirim~\cite{Y2004}, respectively. The optimal partition approach fully describes the optimal set mapping and the optimal value function on the entire $\interior(\mathcal{E})$. In contrast to the optimal basis approach in LO~\cite{JRT93}, which may produce inconsistent results due to problem degeneracy, the results from the optimal partition approach is unique and invariant with respect to any regularity condition for parametric conic optimization problems. Recently, the second and fourth authors~\cite{MT20} expanded on the optimal partition approach and an invariancy interval in~\cite{GS99} by introducing the concepts of a nonlinearity interval and a transition point for the optimal partition of $(\mathrm{P_{\epsilon}})$ and $(\mathrm{D_{\epsilon}})$. An invariancy interval, see Definition~\ref{def:invariancy_interval}, is an open maximal subinterval of $\interior(\mathcal{E})$ on which the optimal partition is invariant with respect to $\epsilon$. A nonlinearity interval, see Definition~\ref{def:nonlinearity_interval}, is an open maximal subinterval of $\interior(\mathcal{E})$ 
on which the rank of maximally complementary optimal solutions $X^*(\epsilon)$ and $S^*(\epsilon)$ stay constant, while the optimal partition varies with $\epsilon$. A transition point, see Definition~\ref{def:transition_point}, is the boundary point of an invariancy or a nonlinearity interval which belongs to $\interior(\mathcal{E})$. Unlike a parametric LO problem~\cite{JRT93}, the optimal value function of SDO consists of nonlinear pieces (of not necessarily polynomial type) on nonlinearity intervals.
\subsection{Contributions}
Very little is known yet about the nonlinearity intervals and the topology of their optimal solutions for a parametric SDO problem. In particular, in contrast to a parametric LO problem, there is no procedure for the full decomposition of $\interior(\mathcal{E})$ into invariancy and nonlinearity intervals. Our main contribution is a numerical algebraic geometry procedure for the computation of nonlinearity intervals and transition points in $\interior(\mathcal{E})$. To the best of our knowledge, this is the first comprehensive methodology for the full decomposition of $\interior(\mathcal{E})$ for a parametric SDO problem. 

\vspace{5px}
\noindent
The first part of this paper reviews the notions of invariancy set, nonlinearity interval, and transition point and investigates their characterizations. We prove that the set of transition points is finite, see Theorem~\ref{transition_finite}, and using continuity arguments on the basis of Painlev\'e-Kuratowski set convergence, we provide sufficient conditions under which a nonlinearity interval exists, see Lemma~\ref{continuity_nonlinearity}. 
We analyze the continuity of the optimal set mapping and show that continuity may fail on a nonlinearity interval, see Example~\ref{elliptope_example}.
Additionally, we show that even a continuous selection~\cite[Chapter 5(J)]{Rock09} through the relative interior of the optimal sets might fail to exist, see problem~\eqref{3elliptope_cut_eq}. The second part of this paper investigates the computation of nonlinearity intervals and transition points of the optimal partition. Under a local nonsingularity condition, see Theorem~\ref{thm:V_ODE}, we develop a methodology, Algorithms~\ref{singular_identification} and~\ref{transition_identification}, to compute the boundary points of a nonlinearity interval and identify a transition point. By assuming a \textcolor{black}{generic} global nonsingularity condition, see Proposition~\ref{extrema2singular}, we then present a numerical procedure, Algorithm~\ref{outermost_alg}, which partitions $\interior(\mathcal{E})$ into a finite union of invariancy intervals, nonlinearity intervals, and transition points. 

\vspace{5px}
\noindent
Since the maximal rank of optimal solutions is preserved on invariancy and nonlinearity intervals, our numerical procedure could be of great interest to the parametric analysis of matrix completion problems, see e.g.,~\cite{AW00}.
Besides sensitivity analysis purposes and their economical interpretations, the identification of a nonlinearity interval is important from practical perspectives. 
For example, in order to approximate the optimal value function 
on a neighborhood of a given $\epsilon$, one needs to 
utilize samples from the same nonlinearity interval containing $\epsilon$. Cifuentes et al.~\cite{CAPT17} studied the local stability of SDO relaxations for polynomial and semi-algebraic optimization problems with emphasis on a notion similar to a nonlinearity interval.

\subsection{Organization of the paper}\label{organization}
The rest of this paper is organized as follows. In Section~\ref{continuity_SDO}, we investigate the continuity of the feasible and optimal set mappings at a given $\epsilon \in \interior(\mathcal{E})$ relative to $\interior(\mathcal{E})$. In Section~\ref{parametric_optimal_partition}, we study the sensitivity of the optimal partition with respect to $\epsilon$. Further, we use continuity and semi-algebraicity arguments to characterize nonlinearity intervals and transition points, and we investigate the continuity of the optimal set mapping on a nonlinearity interval. In Section~\ref{nonlinearity_interval_algorithm}, we present an algorithm to compute invariancy intervals, nonlinearity intervals, and transition points in $\interior(\mathcal{E})$. Our numerical experiments are presented in Section~\ref{numerical_experiments}. Finally, we present remarks and topics for future research in Section~\ref{conclusion}.

\paragraph{Notation} Throughout this paper, $\mathbb{S}^n_+$ denotes the cone of $n \times n$ positive semidefinite matrices, $\bd(\cdot)$ represents the boundary of a set, and $\|\cdot\|_2$ denotes the $\ell_2$ norm of a vector. Associated with a symmetric matrix $X$, $\lambda_{\min}(X)$ denotes the smallest eigenvalue of $X$, $\kernel(X)$ is the null space of $X$, and $\svectorize(X)$ denotes a linear mapping stacking the upper triangular part of a symmetric matrix, in which the off-diagonal entries are multiplied by $\sqrt{2}$, i.e.,
 \begin{align}\label{isomorphism}
 \svectorize(X)\!:=\!\big(X_{11}, \sqrt{2}X_{12},\ldots, \sqrt{2}X_{1n}, X_{22}, \sqrt{2} X_{23},\ldots,\sqrt{2}X_{2n},\ldots, X_{nn}\big)^T.
 \end{align} 
For brevity, we often use the notation $\mathcal{A}\!:=\!\big(\svectorize(A^1),\ldots,\svectorize(A^m)\big)^T$ for a compact representation of the coefficient matrices. Finally, for any two square matrices $K_1$ and $K_2$ and a symmetric matrix $H$, the \textit{symmetric Kronecker product}, denoted by $\otimes_s$, is defined as
 \begin{align*}
 (K_1 \otimes_s K_2)\svectorize(H)\!:= \frac12 \svectorize\big(K_2 H K_1^T + K_1HK_2^T\big),
 \end{align*}
see e.g.,~\cite{Kl02} for more details.

\section{Continuity of the feasible set and optimal set mappings}\label{continuity_SDO}
This section investigates the continuity of the primal and dual feasible set mappings and the outer semicontinuity of the primal and dual optimal set mappings for $(\mathrm{P}_{\epsilon})$ and $(\mathrm{D}_{\epsilon})$. We adopt the notions and definitions from~\cite{RD14,Rock09}.

\vspace{5px}
\noindent
Let $\mathbb{R}^{q}$ and $\mathbb{R}^{l}$ be finite-dimensional Euclidean spaces. A mapping $\Phi:\mathbb{R}^q \rightrightarrows \mathbb{R}^{l}$ is called a \textit{set-valued mapping} if it assigns a subset of $\mathbb{R}^{l}$ to each element of $\mathbb{R}^{q}$. The domain of a set-valued mapping $\Phi$ is $\dom(\Phi)\!:=\!\{\xi: \Phi(\xi) \neq \emptyset \}$, and the range of $\Phi$ is defined as $\mathrm{range}(\Phi)\!:= \!\{\nu: \exists \  \xi \ \st \ \nu \in \Phi(\xi)\}$.

\vspace{5px}
\noindent
The following discussion concisely reviews the continuity of a set-valued mapping on the basis of Painlev\'e-Kuratowski set convergence, see~\cite[Chapters 4 and 5]{Rock09} for more details. For a sequence $\{\mathcal{C}_k\}_{k=1}^{\infty}$ of subsets of $\mathbb{R}^l$, the \textit{outer} and \textit{inner} limits are defined, respectively, as
 \begin{align}
\limsup\limits_{k \to \infty} \mathcal{C}_k&\!:=\!\Big \{\nu : \liminf\limits_{k \to \infty} \ \distance(\nu,\mathcal{C}_k)=0 \Big \},\nonumber\\
\liminf\limits_{k \to \infty} \mathcal{C}_k&\!:=\!\Big \{\nu : \limsup\limits_{k \to \infty} \ \distance(\nu,\mathcal{C}_k)=0 \Big \}, \label{liminf_def}
 \end{align}
where $\distance(\nu,\mathcal{C}_k)=\inf_{x \in \mathcal{C}_k} \|\nu-x\|_2$.
Let $\mathcal{X}$ be a subset of $\mathbb{R}^q$ containing $\bar{\xi}$. A set-valued mapping $\Phi$ is called \textit{outer semicontinuous} at $\bar{\xi}$ relative to $\mathcal{X}$ if $\limsup\limits_{\xi \to \bar{\xi}} \Phi(\xi) \subseteq \Phi(\bar{\xi})$ and \textit{inner semicontinuous} at $\bar{\xi}$ relative to $\mathcal{X}$ if $\liminf\limits_{\xi \to \bar{\xi}} \Phi(\xi) \supseteq \Phi(\bar{\xi})$, where  
\begin{align*}
\limsup_{\xi \to \bar{\xi}} \Phi(\xi)&\!:=\!\bigcup_{ \mathcal{X} \supseteq  \xi_k \to \bar{\xi}} \limsup_{k \to \infty} \Phi(\xi_k),\\
\liminf\limits_{\xi \to \bar{\xi}} \Phi(\xi)&\!:=\!\bigcap_{ \mathcal{X} \supseteq  \xi_k \to \bar{\xi}} \liminf_{k \to \infty} \Phi(\xi_k). 
\end{align*}
When $\mathcal{X} = \mathbb{R}^q$, we simply call $\Phi$ outer or inner semicontinuous at $\bar{\xi}$.
\begin{definition}
A set-valued mapping $\Phi$ is 
\textit{Painlev\'e-Kuratowski continuous} at $\bar{\xi}$ relative to $\mathcal{X}$ if it is both outer and inner semicontinuous at $\bar{\xi}$ relative to $\mathcal{X}$.
\end{definition}

\vspace{5px}
\noindent
In our setting, outer and inner semicontinuity agree with the notions of closedness and openness of a point-to-set map in~\cite{Ho73b}, see also~\cite[Theorem 5.7(c)]{Rock09} and~\cite[Corollary 1.1]{Ho73b}.

\vspace{5px}
\noindent
We show the continuity of the feasible set mapping and the
outer semicontinuity of the optimal set mapping relative to $\interior(\mathcal{E})$. 
Trivially, $\mathcal{P}:\mathbb{R} \rightrightarrows \mathbb{S}^n$ is
continuous since it remains invariant with respect to $\epsilon$.
Furthermore, the continuity of $\mathcal{D}:\mathbb{R} \rightrightarrows \mathbb{R}^m \times \mathbb{S}^n$ relative to $\interior(\mathcal{E})$ follows from~\cite[Theorems 10 and 12]{Ho73b}, where $\mathcal{D}(\epsilon) = \emptyset$ for every $\epsilon \in \mathbb{R} \setminus \mathcal{E}$, see also~\cite[Example~5.10]{Rock09}. For the sake of completeness, we provide a proof for our special case here. 
\begin{proposition}\label{set_valued_inner}
Under Assumption~\ref{IPC}, the set-valued mapping $\mathcal{D}$ is continuous relative to $\interior(\mathcal{E})$. 
\end{proposition}
\proof{Proof.}
For the sake of brevity, we define $L(y):=\sum_{i=1}^m y_iA^i$. The outer semicontinuity of $\mathcal{D}$ is immediate from the closedness of $\mathbb{S}^n_+$, see e.g.,~\cite[Example 5.8]{Rock09}. Hence, it only remains to show that $\mathcal{D}$ is inner semicontinuous at every $\epsilon' \in \interior(\mathcal{E})$, i.e., given a sequence $\{\epsilon_k\}_{k=1}^\infty$ with 
$\epsilon_k \to \epsilon'$ and an arbitrary $(\hat{y},\hat{S}) \in \mathcal{D}(\epsilon')$, there exists a convergent sequence $(y_k,S_k) \to (\hat{y},\hat{S})$ such that $(y_k,S_k) \in \mathcal{D}(\epsilon_k)$ for all sufficiently large $k$. To that end, let us define $y_k:=(1-\alpha_k)\hat{y} + \alpha_k \bar{y}$ and $S_k\!:=C+\epsilon_k \bar{C}-L(y_k)$, where $(\bar{y},\bar{S}) \in \mathcal{D}(\epsilon')$ such that $\bar{S} \succ 0$. By Assumption~\ref{IPC}, such a $(\bar{y},\bar{S})$ exists. We then need to construct a convergent sequence $\alpha_k \to 0$ such that $S_k \succeq 0$ holds. We assume that $\lambda_{\min}(\hat{S}) = 0$, since otherwise for any arbitrary sequence $\alpha_k \to 0$ we always have $S_k \succ 0$ when $k$ is sufficiently large. 

\vspace{5px}
\noindent
Notice that if $0 \le \alpha_k \le 1$, then $S_k \succeq 0$ is satisfied by requiring
\begin{align*}
(1-\alpha_k)\lambda_{\min} \big (C+\epsilon_k \bar{C}- L(\hat{y}) \big)+\alpha_k \lambda_{\min} \big (C+\epsilon_k \bar{C}-L(\bar{y}) \big ) \ge 0,
\end{align*}
which is equivalent to
\begin{align*}
\alpha_k \ge \mu_k \!:=\frac{-\lambda_{\min}\big (C+\epsilon_k \bar{C}-L(\hat{y}) \big )}{\lambda_{\min}\big (C+\epsilon_k \bar{C}-L(\bar{y}) \big)-\lambda_{\min} \big (C+\epsilon_k \bar{C}-L(\hat{y}) \big )}
\end{align*}
for sufficiently large $k$, since the denominator has to be positive. Letting $\alpha_k\!:=\max\{\mu_k,0\}$, we get the desired sequence.  \Halmos
\endproof

\vspace{5px}
\noindent
As a result of Proposition~\ref{set_valued_inner}, we can show that
$\mathcal{P}^*:\mathbb{R} \rightrightarrows \mathbb{S}^n$ and $\mathcal{D}^*:\mathbb{R} \rightrightarrows \mathbb{R}^m \times \mathbb{S}^n$ are outer semicontinuous relative to $\interior(\mathcal{E})$, see e.g.,~\cite[Theorem 8]{Ho73b} or~\cite[Theorem~3B.5]{RD14}. All this implies that for any $\epsilon' \in \interior(\mathcal{E})$ and any sequence $\epsilon_k \to \epsilon'$ we have
\begin{align}\label{optimal_set_semicontinuity}
\liminf_{k \to \infty} \mathcal{P}^*(\epsilon_k) \subseteq \limsup_{k \to \infty} \mathcal{P}^*(\epsilon_k) \subseteq \mathcal{P}^*(\epsilon') \ \ \ \text{and} \ \ \
\liminf_{k \to \infty} \mathcal{D}^*(\epsilon_k) \subseteq \limsup_{k \to \infty} \mathcal{D}^*(\epsilon_k) \subseteq \mathcal{D}^*(\epsilon').
\end{align}
However, $\mathcal{P}^*$ and $\mathcal{D}^*$ are not necessarily inner semicontinuous relative to $\interior(\mathcal{E})$ as shown
in Example~\ref{elliptope_example}, where the optimal set is multiple-valued at $\epsilon=\frac12$ but single-valued everywhere else in a neighborhood of $\tfrac12$. 
Nevertheless, the set of points at which $\mathcal{P}^*$ or $\mathcal{D}^*$ fails to be continuous relative to $\interior(\mathcal{E})$ is of \textit{first category} in $\interior(\mathcal{E})$, i.e., it is the union of countably many nowhere dense sets in $\interior(\mathcal{E})$, see e.g.,~\cite{M2000}. This directly follows from the outer semicontinuity of the optimal set mapping relative to $\interior(\mathcal{E})$ and Theorem 5.55 in~\cite{Rock09}. All this yields the following result.
\begin{proposition}\label{number_of_discontinuous_points}
The set of points at which $\mathcal{P}^*$ or $\mathcal{D}^*$ fails to be continuous relative to $\interior(\mathcal{E})$ has empty interior. 
\end{proposition}
\proof{Proof.}
Since $\interior(\mathcal{E})$ is a Baire subset of $\mathbb{R}$~\cite[Lemma 48.4]{M2000}, every first category subset of $\interior(\mathcal{E})$ has empty interior. \Halmos
\endproof
\noindent
As a consequence of Proposition~\ref{number_of_discontinuous_points}, every open subset of $\interior(\mathcal{E})$ contains a point at which both $\mathcal{P}^*$ and $\mathcal{D}^*$ are continuous relative to $\interior(\mathcal{E})$.
\section{Sensitivity of the optimal partition}\label{parametric_optimal_partition}
We briefly review the notions of an invariancy interval, nonlinearity interval, and a transition point from~\cite{MT20}. Let $\pi(\epsilon)\!:=\!\big(\mathcal{B}(\epsilon),\mathcal{T}(\epsilon),\mathcal{N}(\epsilon)\big)$ denote the subspaces of the optimal partition at $\epsilon$, and let $\big(Q_{\mathcal{B}(\epsilon)},Q_{\mathcal{T}(\epsilon)}, Q_{\mathcal{N}(\epsilon)}\big)$ be an orthonormal basis partitioned according to the subspaces of the optimal partition.

\begin{definition}[\cite{GS99,MT20}]\label{def:invariancy_interval}
An \textit{invariancy set} is a maximal subset $\mathcal{I}_{\mathrm{inv}}$ of $\interior(\mathcal{E})$ on which $\pi(\epsilon)$ is invariant for all $\epsilon \in \mathcal{I}_{\mathrm{inv}}$. 
\end{definition}
Indeed, an invariancy set is proven to be either a singleton or an open, possibly unbounded, subinterval of $\interior(\mathcal{E})$, see~\cite[Lemma~3.3]{MT20} and its preceding discussion. A non-singleton $\mathcal{I}_{\mathrm{inv}}$ is simply called an \textit{invariancy interval}. 
\begin{remark}
Even though the optimal partition of a singleton $\mathcal{I}_{\mathrm{inv}}$ is vacuously invariant on $\mathcal{I}_{\mathrm{inv}}$, it differs from the optimal partition of every neighborhood of $\mathcal{I}_{\mathrm{inv}}$. \halmos
\end{remark}
\noindent
The primal optimal set mapping $\mathcal{P}^*$ is constant on an invariancy interval~\cite[Remark~3.1]{MT20}. Furthermore, the boundary points of an invariancy set, containing a given $\bar{\epsilon}$, can be efficiently computed by solving a pair of auxiliary SDO problems~\cite[Lemma 4.1]{GS99}:
\begin{equation}\label{auxiliary_problems}
\begin{aligned}
\alpha_{\mathrm{inv}}(\beta_{\mathrm{inv}})\!:=\inf(\sup)  &\quad \epsilon\\[-1\jot] 
\st &\quad \sum_{i=1}^m y_i A^i + Q_{\mathcal{N}(\bar{\epsilon})} U_S Q^T_{\mathcal{N}(\bar{\epsilon})} = C + \epsilon \bar{C},\\
&\quad U_S \succeq 0,
\end{aligned}
\end{equation}

\noindent
where we might have $\alpha_{\mathrm{inv}} = -\infty$, $\beta_{\mathrm{inv}} = \infty$, or both. If $\alpha_{\mathrm{inv}} < \bar{\epsilon} < \beta_{\mathrm{inv}}$ holds, then $\bar{\epsilon}$ belongs to an invariancy interval. Otherwise, $\bar{\epsilon}$ belongs to a nonlinearity interval, or it is a transition point, as formally defined in Definitions~\ref{def:nonlinearity_interval} and~\ref{def:transition_point}. Recall that $\big(X^*(\epsilon),y^*(\epsilon),S^*(\epsilon)\big)$ denotes a maximally complementary optimal solution.

\begin{definition}[Definition 3.6 in~\cite{MT20}]\label{def:nonlinearity_interval}
A \textit{nonlinearity interval} is an open maximal subinterval $\mathcal{I}_{\mathrm{non}}$ of $\interior(\mathcal{E})$ 
on which both $\rank\!\big(X^*(\epsilon)\big)$ and $\rank\!\big(S^*(\epsilon)\big)$ are constant while $\pi(\epsilon)$ varies with $\epsilon$, i.e., $\epsilon_1 \neq \epsilon_2$ implies $\pi(\epsilon_1) \neq \pi(\epsilon_2)$ for all $\epsilon_1,\epsilon_2 \in \mathcal{I}_{\mathrm{non}}$.
\end{definition}

\begin{definition}[Definition 3.5 in~\cite{MT20}]\label{def:transition_point}
A point $\bar{\epsilon} \in \interior(\mathcal{E})$ is called a \textit{transition point} if for every $\delta > 0$, there exists $\epsilon \in (\bar{\epsilon} - \delta, \bar{\epsilon} + \delta) \cap \interior(\mathcal{E})$ such that
\begin{align*}
\rank\!\big(X^*(\epsilon)\big)\neq \rank\!\big(X^*(\bar{\epsilon})\big) \quad \text{or} \quad \rank\!\big(S^*(\epsilon)\big)\neq \rank\!\big(S^*(\bar{\epsilon})\big).
\end{align*} 
\end{definition}
Definition~\ref{def:transition_point} is consistent with the one defined for a parametric LO problem~\cite{JRT93}, as spelled out in the following proposition. 
\begin{proposition}\label{transition_point_equivalency}
At a boundary point $\bar{\epsilon} \in \interior(\mathcal{E})$ of an invariancy interval $\mathcal{I}_{\mathrm{inv}}$ and for some $\hat{\epsilon} \in \mathcal{I}_{\mathrm{inv}}$ we have
\begin{align*}
\rank\!\big(X^*(\hat{\epsilon})\big)\neq \rank\!\big(X^*(\bar{\epsilon})\big) \quad \text{or} \quad \rank\!\big(S^*(\hat{\epsilon})\big)\neq \rank\!\big(S^*(\bar{\epsilon})\big).
\end{align*}
\end{proposition}
Before proving this statement, we need the following result. 
\begin{proposition}\label{primal_constancy}
If $\mathcal{P}^*$ and $\rank\big(S^*(\epsilon)\big)$ are constant on $[\epsilon_1,\epsilon_2]$, then so is $\pi(\epsilon)$.
\end{proposition}
\proof{Proof.}
Let us define $\epsilon_{\gamma}\!:=\!\gamma \epsilon_1 + (1-\gamma) \epsilon_2$, where $\gamma \in [0,1]$. Then for every $\gamma \in (0,1)$ it is easy to verify that $\big(X(\epsilon_{\gamma}),y(\epsilon_{\gamma}),S(\epsilon_{\gamma})\big)$ is an optimal solution of $(\mathrm{P_{\epsilon_{\gamma}}})-(\mathrm{D_{\epsilon_{\gamma}}})$, where
\begin{align}\label{intermediate_maximally_complementary}
X(\epsilon_{\gamma})\!:= X^*(\epsilon_1), \quad y(\epsilon_{\gamma})\!:= \gamma y^*(\epsilon_1) + (1-\gamma) y^*(\epsilon_2),\quad S(\epsilon_{\gamma})\!:= \gamma S^*(\epsilon_1) + (1-\gamma) S^*(\epsilon_2),
\end{align}
in which $X(\epsilon_{\gamma})S(\epsilon_{\gamma})=0$ follows from the constancy of $\mathcal{P}^*$. Let $0 < \gamma_1,\gamma_2 < 1$. Notice from~\eqref{intermediate_maximally_complementary} and from the positive semidefiniteness of $S^*(\epsilon_1)$ and $S^*(\epsilon_2)$ that for every $q \in \mathbb{R}^n$, $q^T S(\epsilon_{\gamma_1})q = 0$ implies 
\begin{align*}
q^T S^*(\epsilon_{1})q = 0 \quad \text{and} \quad q^T S^*(\epsilon_{2})q = 0, 
\end{align*}
which in turn yield $q^T S(\epsilon_{\gamma_{2}})q = 0$ by~\eqref{intermediate_maximally_complementary}. Therefore, $\kernel\!\big(S(\epsilon_{\gamma_1})\big) \subseteq \kernel\!\big(S(\epsilon_{\gamma_2})\big)$, and by switching the roles of $\gamma_1$ and $\gamma_2$ we get $\kernel\!\big(S(\epsilon_{\gamma_1})\big) = \kernel\!\big(S(\epsilon_{\gamma_2})\big)$. Further, it is obvious from~\eqref{intermediate_maximally_complementary} that $\kernel\!\big(X(\epsilon_{\gamma_{1}})\big) = \kernel\!\big(X(\epsilon_{\gamma_2})\big)$. Finally, we can conclude from the constancy of the primal optimal set and $\rank\!\big(S^*(\epsilon)\big)$ on $[\epsilon_1,\epsilon_2]$ that $\rank\!\big(X(\epsilon_{\gamma})\big) = \rank\!\big(X^*(\epsilon_{\gamma})\big)$ and $\rank\!\big(S(\epsilon_{\gamma})\big) \ge \rank\!\big(S^*(\epsilon_{\gamma})\big)$ for all $\gamma \in (0,1)$, which in turn indicate that $\big(X(\epsilon_{\gamma}),y(\epsilon_{\gamma}),S(\epsilon_{\gamma})\big)$ is maximally complementary. \halmos
\endproof

\proof{Proof of Proposition~\ref{transition_point_equivalency}.}
In addition to Proposition~\ref{primal_constancy}, we need to recall from~\eqref{optimal_set_semicontinuity} that for any sequence $\mathcal{I}_{\mathrm{inv}} \supseteq \epsilon_k \to \bar{\epsilon}$, it holds that $\liminf_{k \to \infty} \mathcal{P}^*(\epsilon_k) \subseteq \mathcal{P}^*(\bar{\epsilon})$, whereas $\liminf_{k \to \infty} \mathcal{P}^*(\epsilon_k) = \mathcal{P}^*(\hat{\epsilon})$ follows from the constancy of $\mathcal{P}^*$ on $\mathcal{I}_{\mathrm{inv}}$ and~\cite[Exercise 4.3(b)]{Rock09}. Consequently, $\mathcal{P}^*(\hat{\epsilon}) \subseteq \mathcal{P}^*(\bar{\epsilon})$, and exactly one of the following holds: (a) $ \mathcal{P}^*(\hat{\epsilon}) \subseteq \bd\!\big(\mathcal{P}^*(\bar{\epsilon})\big)$ or (b) $\mathcal{P}^*(\hat{\epsilon}) \cap \ri\!\big(\mathcal{P}^*(\bar{\epsilon})\big) \neq \emptyset$. Case (a) leads to $\rank\!\big(X^*(\hat{\epsilon})\big) < \rank\!\big(X^*(\bar{\epsilon})\big)$ by the definition of a maximally complementary optimal solution, while case (b) implies $\ri\!\big(\mathcal{P}^*(\hat{\epsilon})\big) \subseteq \ri\!\big(\mathcal{P}^*(\bar{\epsilon})\big)$ and thus $\rank\!\big(S^*(\hat{\epsilon})\big)\neq \rank\!\big(S^*(\bar{\epsilon})\big)$ by the proof of Proposition~\ref{primal_constancy}. \halmos
\endproof
\begin{remark}\label{optimal_set_behavior}
It is immediate from Proposition~\ref{primal_constancy} that on a nonlinearity interval both the primal and dual optimal sets must vary with $\epsilon$.  \halmos
\end{remark}
\noindent
A boundary point of an invariancy or a nonlinearity interval, if it belongs to $\interior(\mathcal{E})$, must be a transition point by Definition~\ref{def:nonlinearity_interval} and Proposition~\ref{transition_point_equivalency}. On the other hand, the semi-algebraic~\cite{BPR06} property of Definitions~\ref{def:invariancy_interval} and~\ref{def:nonlinearity_interval} implies that the set of transition points is always finite, see Theorem~\ref{transition_finite}, i.e., a transition point must be a boundary point of an invariancy or a nonlinearity interval. The idea of the proof is analogous to~\cite[Theorem~1]{MT21} for the optimal partition of a parametric second-order conic optimization problem. For the sake of completeness, we refer the reader to the Appendix for a self-contained proof.

\begin{theorem}\label{transition_finite}
The set of transition points is finite.
\end{theorem}
\noindent
As a result of Theorem~\ref{transition_finite}, $\interior(\mathcal{E})$ can be always partitioned into the finite union of invariancy intervals, nonlinearity intervals, and transition points. The following example is adopted from~\cite[Example 3.1]{MT20} and shows the existence of nonlinearity intervals and transition points.
\begin{example}\label{elliptope_example}
Consider the following parametric convex optimization problem:
\begin{align}\label{3elliptope}
\min\Bigg\{ (4\epsilon-2)x+(2-4\epsilon)y-2z : \begin{pmatrix} 1 & \ x & \ y\\x & \ 1 & \ z\\y & \ z & \ 1\end{pmatrix} \succeq 0\Bigg\},
\end{align} 
in which the feasible region is a 3-elliptope~\cite{BPT13}, see Figure~\ref{fig:3elliptope}. Since the perturbation parameter $\epsilon$ only appears in the objective function, we can cast the parametric problem~\eqref{3elliptope} into the primal form $(\mathrm{P_{\epsilon}})$ with $X \in \mathbb{S}^3$ and $m=3$ by introducing
\begin{align*}
A^1&=\begin{pmatrix} 1 & 0 & 0\\0 & 0  & 0\\0 & 0 & 0\end{pmatrix}, & A^2&=\begin{pmatrix} 0 & 0 & 0\\0 & 1  & 0\\0 & 0 & 0\end{pmatrix}, & A^3&=\begin{pmatrix} 0 & 0 & 0\\0 & 0  & 0\\0 & 0 & 1\end{pmatrix}, \\
C&=\begin{pmatrix} \ \ 0 & -1 & \ \ 1\\ -1 & \ \ 0  & -1\\ \ \ 1 & -1 & \ \ 0\end{pmatrix}, & \bar{C}&=\begin{pmatrix} \ \ 0 & \ 2 & -2\\ \ \ 2 & \ 0  & \ \ 0\\-2 & \ 0 & \ \ 0\end{pmatrix}, & b&=(1, \ 1, \ 1)^T.
\end{align*}
For all $\epsilon \in (-\tfrac12,\tfrac32)$, see~\cite[Example 3.1]{MT20}, a strictly complementary optimal solution is given by
\begin{figure}
 \begin{minipage}[c]{0.5\textwidth}
\includegraphics[height=1.8in]{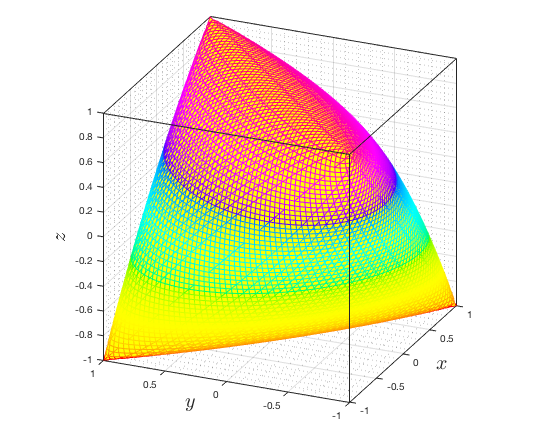}
\end{minipage}\hfill
\begin{minipage}[c]{0.5\textwidth}
\caption{The feasible set of the parametric convex optimization problem~\eqref{3elliptope}, being invariant with respect to $\epsilon$.}
\label{fig:3elliptope}
\end{minipage}
\end{figure}
\small
\begin{align*}
X^*(\epsilon)=\begin{pmatrix} 1 & \frac12 - \epsilon & \epsilon - \frac12 \\ \frac12-\epsilon &  1 &  1-2(\epsilon-\frac12)^2\\ \epsilon-\frac12 & 1-2(\epsilon-\frac12)^2 & 1 \end{pmatrix}, \ y^*(\epsilon)=\begin{pmatrix} -(2\epsilon-1)^2\\ -1\\ -1 \end{pmatrix}, \ S^*(\epsilon)=\begin{pmatrix} (2\epsilon-1)^2 &  2\epsilon-1 & 1- 2\epsilon \\ 2\epsilon-1 & \ \ 1 & -1\\ 1- 2\epsilon & -1 & \ \ 1 \end{pmatrix},
\end{align*}
\normalsize
while a maximally complementary optimal solution at $\epsilon=\frac32$ is given by
\begin{align*}
X^*(\tfrac32)&=\begin{pmatrix} \ \ 1 & -1 & \ \ 1 \\ -1 &  \ \ 1 &  -1\\ \ \ 1 & -1 & \ \ 1  \end{pmatrix}, \ y^*(\tfrac32)=\begin{pmatrix} -4 \\ -1\\-1 \end{pmatrix}, \ S^*(\tfrac32)=\begin{pmatrix} \ \ 4 & \ \ 2 & -2 \\ \ \ 2 & \ \ 1 & -1 \\ -2 & -1 & \ \ 1 \end{pmatrix}.
\end{align*}
The eigenvalue decompositions of $X^*(\epsilon)$ and $S^*(\epsilon)$ reveal that
\begin{align*}
\rank\!\big(X^*(\epsilon)\big)=\begin{cases} 2 \qquad \epsilon \in (-\frac12,\frac32),\\1 \qquad \epsilon = \frac32, \end{cases} \qquad \rank\!\big(S^*(\epsilon)\big)= 1, \qquad \epsilon \in (-\tfrac12,\tfrac32]. 
\end{align*}
By definition, $(-\frac12,\frac32)$ is a nonlinearity interval and $\epsilon = \frac32$ is a transition point of the optimal partition. \halmos
\end{example}

\vspace{5px}
\noindent
Due to unknown behavior of the optimal set mapping in a parametric SDO problem, see Remark~\ref{optimal_set_behavior}, a general existence condition for a nonlinearity interval or a transition point is still an open question. Nevertheless, strict complementarity coupled with the continuity of the optimal set mapping at a given $\bar{\epsilon}$ relative to $\interior(\mathcal{E})$
provide sufficient conditions for the existence of a 
nonlinearity interval surrounding~$\bar{\epsilon}$.
\begin{lemma}\label{continuity_nonlinearity}
Let $\{\bar{\epsilon}\}$ be a singleton invariancy set, and let $\big(X^*(\bar{\epsilon}),y^*(\bar{\epsilon}),S^*(\bar{\epsilon})\big)$ be a strictly complementary optimal solution at $\bar{\epsilon} \in \interior(\mathcal{E})$, at which both the primal and dual optimal set mappings are continuous relative to $\interior(\mathcal{E})$.
Then $\bar{\epsilon}$ belongs to a nonlinearity interval.
\end{lemma}
\proof{Proof.}
The strict complementarity condition yields
\begin{align*}
\rank\!\big(X^*(\bar{\epsilon})\big) + \rank\!\big(S^*(\bar{\epsilon})\big) = n.
\end{align*}
Continuity of $\mathcal{P}^*$ and $\mathcal{D}^*$ at $\bar{\epsilon}$, along with the continuity of the eigenvalues, shows that $\rank\!\big(X^*(\bar{\epsilon})\big)  \le \rank\!\big(X^*(\epsilon)\big)$ and $\rank\!\big(S^*(\bar{\epsilon})\big) \le \rank\!\big(S^*(\epsilon)\big)$ for all $\epsilon$ in a small neighborhood 
of $\bar{\epsilon}$, see also~\cite[Theorem 3B.2(b)]{RD14}. Hence, the rank of $X^*(\epsilon)$ and $S^*(\epsilon)$ remain constant on a sufficiently small neighborhood of $\bar{\epsilon}$.  \Halmos
\endproof

\vspace{5px}
\noindent
Unfortunately, the converse of Lemma~\ref{continuity_nonlinearity} is not necessarily true. In fact, the primal or dual optimal set mapping might fail to be continuous on a nonlinearity interval. This can occur since the $\liminf\limits$ of a sequence of faces is not necessarily a face of the feasible set, i.e., it might be a subset of the relative interior of a face. A counterexample is Example~\ref{elliptope_example}, where the strict complementarity condition holds on a nonlinearity interval $(-\frac12,\frac32)$. The primal optimal set mapping is single-valued everywhere on $(-\frac12,\frac12)\cup(\frac12,\frac32)$, see~\cite[Page~204]{MT20}. However, $\mathcal{P}^*$ fails to be inner semicontinuous at $\epsilon=\frac12$, because $\mathcal{P}^*$ is multiple-valued at $\epsilon=\frac12$, and 
\begin{align*}
\liminf\limits_{k \to \infty} \mathcal{P}^*(\epsilon_k) \subset \ri(\mathcal{P}^*(\tfrac12))
\end{align*}
for any sequence $\epsilon_k \to \frac12$. 

\begin{figure}
 \begin{minipage}[c]{0.50\textwidth}
\includegraphics[height=2.0in]{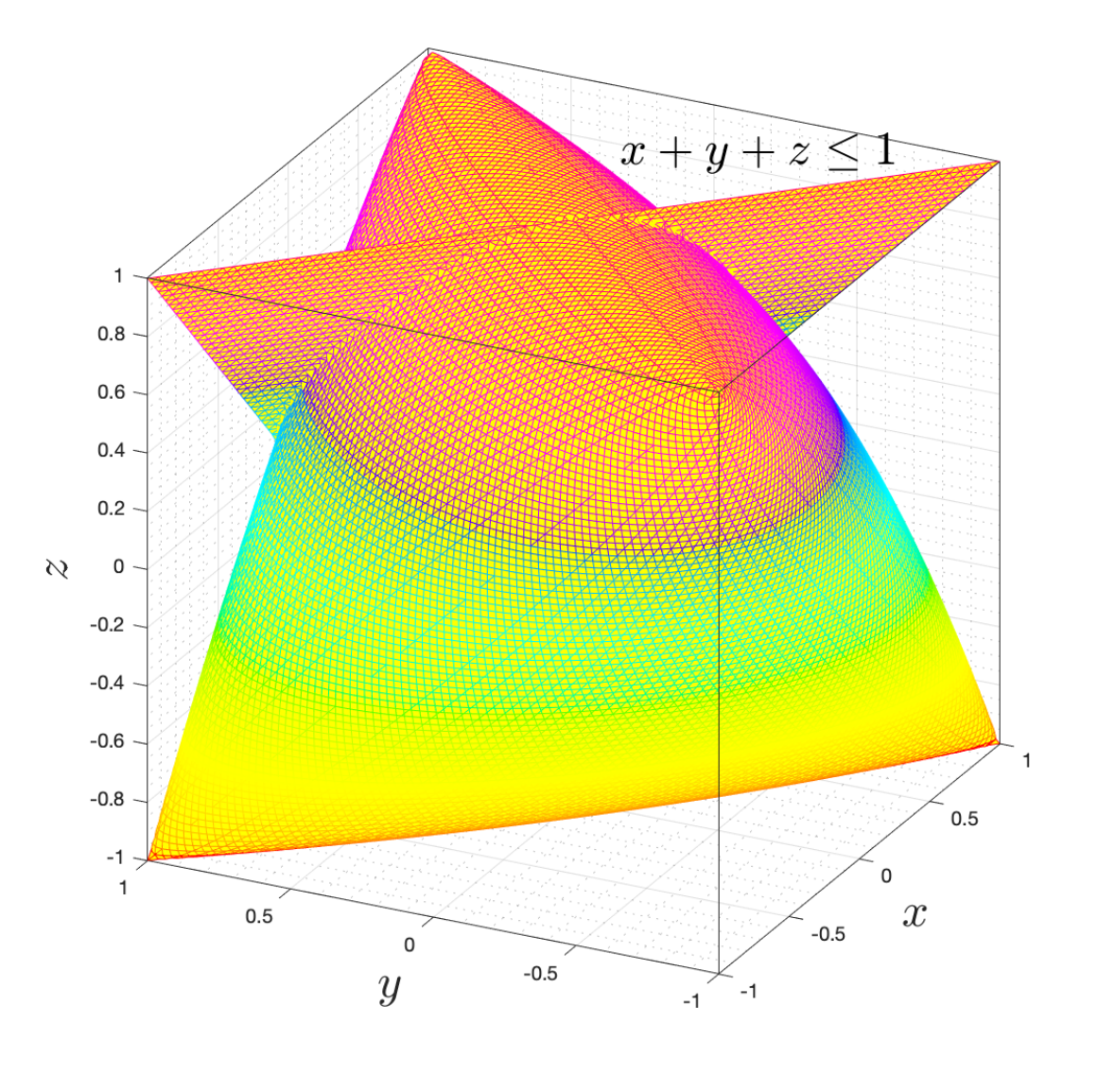}
\end{minipage}\hfill
\begin{minipage}[c]{0.5\textwidth}
\caption{The feasible set of the parametric convex optimization problem~\eqref{3elliptope_cut_eq}.}
\label{fig:3elliptope_cut}
\end{minipage}
\end{figure}

\begin{remark}
The continuity condition in Lemma~\ref{continuity_nonlinearity} can be relaxed by imposing the conditions
\begin{align}\label{continuity_selection}
\liminf_{k \to \infty} \mathcal{P}^*(\epsilon_k) \cap \ri\big(\mathcal{P}^*(\bar{\epsilon})\big) \neq \emptyset\quad \text{and} \quad \liminf_{k \to \infty} \mathcal{D}^*(\epsilon_k) \cap \ri\big(\mathcal{D}^*(\bar{\epsilon})\big) \neq \emptyset
\end{align}
for every sequence $\epsilon_k \to \bar{\epsilon}$, which by~\eqref{liminf_def} and the continuity of the eigenvalues imply the existence of a nonlinearity interval around $\bar{\epsilon}$, see also~\cite[Theorem~3.7]{MT20}. However, even the weaker condition~\eqref{continuity_selection} may not hold on a nonlinearity interval. For instance, by adding the inequality constraint \hbox{$x+y+z \le 1$} to problem~\eqref{3elliptope} we get
\begin{align}\label{3elliptope_cut_eq}
\min\Bigg\{ (4\epsilon-2)x+(2-4\epsilon)y-2z : \begin{pmatrix} 1 & \ x & \ y \\x & \ 1 & \ z \\y & \ z & \ 1\end{pmatrix} \succeq 0, \quad x+y+z \le 1\Bigg\},
\end{align}  
which can be analogously cast into the primal form $(\mathrm{P_{\epsilon}})$ with $X \in \mathbb{S}^4$ and $m=7$, see Figure~\ref{fig:3elliptope_cut}. For all $\epsilon \in (-\frac12,\frac32)\setminus\{\frac12\}$ we still have a unique strictly complementary optimal solution 
\begin{align*}
X^*(\epsilon)&=\begin{pmatrix} 1 & \frac12 - \epsilon & \epsilon - \frac12 & 0 \\ \frac12-\epsilon &  1 &  1-2(\epsilon-\frac12)^2 & 0\\ \epsilon-\frac12 & 1-2(\epsilon-\frac12)^2 & 1 & 0\\0 & 0 & 0 & 2(\epsilon-\frac12)^2 \end{pmatrix}, \ y^*(\epsilon)=\begin{pmatrix} -(2\epsilon-1)^2, \ -1, \ -1, \ 0, \ 0, \ 0, \ 0 \end{pmatrix}^T, \\ S^*(\epsilon)&=\begin{pmatrix} (2\epsilon-1)^2 &  2\epsilon-1 & 1- 2\epsilon & 0 \\ 2\epsilon-1 & \ \ 1 & -1 & 0\\ 1- 2\epsilon & -1 & \ \ 1 & 0\\0 & \ \ 0 & \ \ 0 & 0\end{pmatrix}.
\end{align*}
\noindent
However, for any $\epsilon_k \to \frac12$ the sequence $X^*(\epsilon_k)$ converges to an optimal solution on the boundary of $\mathcal{P}^*(\frac12)$. This example shows that even a continuous selection~\cite[Chapter 5(J)]{Rock09} through the relative interior of the optimal sets might fail to exist on a nonlinearity interval. However, we do not know yet whether~\eqref{continuity_selection} could fail at a boundary point of a nonlinearity interval.  \halmos
\end{remark}

\section{Identification of the optimal partitions}\label{nonlinearity_interval_algorithm}
This section proposes a methodology to compute the boundary points of nonlinearity intervals and identify transition points in $\interior(\mathcal{E})$.
By Theorem~\ref{transition_finite}, the interval $\interior(\mathcal{E})$ is the disjoint union of finitely many invariancy intervals, nonlinearity intervals, and transition points. An invariancy interval can be efficiently computed by solving the auxiliary SDO problems~\eqref{auxiliary_problems}. In general, however, the identification of a nonlinearity interval around a given $\bar{\epsilon}$ is a nontrivial computational task, since the conditions of Lemma~\ref{continuity_nonlinearity} may not be easily checked in practice. One could try to simply solve $(\mathrm{P_{\epsilon}})$ and $(\mathrm{D_{\epsilon}})$ for various $\epsilon$ in a neighborhood of $\bar{\epsilon}$ with the aim of finding the desired nonlinearity interval. However, this approach could fail due to the fact that the solution of IPMs usually come with numerical inaccuracy. Therefore, a positive eigenvalue of $X^*(\epsilon)$ or $S^*(\epsilon)$, which could be doubly exponentially small~\cite[Example~3.2]{MT19}, may not be identified. On the other hand, since the set of transition points is finite, see Theorem~\ref{transition_finite}, the numerical inaccuracy could lead one to miss a transition point when simply solving $(\mathrm{P_{\epsilon}})$ and $(\mathrm{D_{\epsilon}})$ at a given set of~mesh points.

\vspace{5px}
\noindent
In order to compute the boundary points of nonlinearity intervals, we numerically locate the transition points by reformulating the optimality conditions~\eqref{KKT_param_SDO} as a system of polynomials. We then view the problem of finding transition points through the lens of numerical algebraic geometry, see~\cite{BSHW2013,SW05} for an overview of results regarding polynomial systems. 

\subsection{Algebraic formulation}\label{sec:algebraic_formulation}
For $\mathcal{A}\!:=\!\big(\svectorize(A^1),\ldots,\svectorize(A^m)\big)^T$,
the optimality conditions~\eqref{KKT_param_SDO} can be equivalently written as 
\begin{align}\label{KKT_conditions}
F(V,\epsilon)\!:=\!\begin{pmatrix}
\mathcal{A}\svectorize(X) - b\\
\mathcal{A}^Ty + \svectorize(S) - \svectorize(C+\epsilon\bar{C})\\
\frac12 \svectorize(XS + SX)\\
\end{pmatrix}=0,
\end{align}
\begin{align}\label{eq:KKT_INEQ}
X,S \succeq 0,
\end{align}
where $V\!:=\!\big(\svectorize(X);y;\svectorize(S)\big)$ is the vector of variables. Given a particular $\epsilon$, the algebraic set of solutions satisfying~\eqref{KKT_conditions} is denoted by
\begin{align}\label{algebraic_set_optimal}
\mathbf{V}\big(F(V,\epsilon)\big)\!:=\!\big\{V\in\mathbb{C}^{m+2t(n)} : F(V,\epsilon)=0\big\},
\end{align}
where $t(n):=n(n+1)/2$. An algebraic set is the solution set of a system of polynomials over~$\mathbb{C}$. Following this notation, a solution in $\mathbf{V}\big(F(V,\epsilon)\big)$, an optimal solution, and a maximally complementary optimal solution of $(\mathrm{P_{\epsilon}})$ and $(\mathrm{D_{\epsilon}})$ are denoted by $\underline{V}(\epsilon)$, $V(\epsilon)$, and $V^*(\epsilon)$, respectively. Clearly, $\underline{V}(\epsilon)$ is not necessarily an optimal solution of $(\mathrm{P_{\epsilon}})$ and $(\mathrm{D_{\epsilon}})$ since it may be complex or fail to satisfy~\eqref{eq:KKT_INEQ}.

\vspace{5px}
\noindent
The Jacobian matrix of~\eqref{KKT_conditions} is given by
\begin{align*}
J(V,\epsilon)\!:=\!\begin{pmatrix} \mathcal{A} & 0 & 0\\0 & \mathcal{A}^T & I_{t(n)} \\ S \otimes_s I_n & 0 & X \otimes_s I_n \end{pmatrix},
\end{align*} 
where the symmetric Kronecker product $\otimes_s$ is defined in Section~\ref{organization}. If the Jacobian is nonsingular at $(V^*(\bar{\epsilon}),\bar{\epsilon})$, then $V^*(\bar{\epsilon})$ is the unique, non-degenerate~\cite[Definitions~5 and 8]{AHO97}, and strictly complementary optimal solution of $(\mathrm{P_{\bar{\epsilon}}})$ and~$(\mathrm{D_{\bar{\epsilon}}})$.

\begin{lemma}[Theorem~3.1 of~\cite{AHO98} and~\cite{Ha98}]\label{nonsingularity_conditions}
The Jacobian $J\big(V^*(\bar{\epsilon}),\bar{\epsilon}\big)$ is nonsingular if and only if the optimal solution $V^*(\bar{\epsilon})$ is non-degenerate and strictly complementary. 
\end{lemma}
\begin{remark}\label{genericity_local_assumption}
We would like to note that non-degeneracy and strict complementarity at fixed $\epsilon$ and $\bar{C}$ are both generic properties~\cite[Theorems~14 and~15]{AHO97}. Therefore, the existence of a unique optimal solution with a nonsingular Jacobian is also a generic property. \halmos
 \end{remark}

\vspace{5px}
\noindent
When the Jacobian is nonsingular, then the implicit function theorem~\cite[Theorem~10.2.1]{D60} and Lemma~\ref{continuity_nonlinearity}  
describe the continuous behavior of $V^*(\epsilon)$ in a neighborhood of $\bar{\epsilon}$ and induce the existence of an invariancy or a nonlinearity interval around $\bar{\epsilon}$. Consequently, transition points and the points at which $\mathcal{P}^*$ or $\mathcal{D}^*$ fails to be continuous relative to $\interior(\mathcal{E})$ are both subsets
 of \textit{singular points} for polynomial system~\eqref{KKT_conditions}, i.e., the set of points
\begin{align*}
\Big\{\epsilon\in\mathbb{C} : \exists  \ \underline{V}(\epsilon) \in \mathbf{V}\big(F(V,\epsilon)\big) \ \text{where the matrix} \ J\big(\underline{V}(\epsilon),\epsilon\big) \ \text{is singular}\Big\},
\end{align*}

\noindent
in which case $\underline{V}(\epsilon)$ is called a \textit{singular solution}. This inclusion might be strict as demonstrated by Example~\ref{elliptope_example}, where $\epsilon=\frac12$
is a singular non-transition point. If $\epsilon$ is not a singular point, then it is called a \textit{nonsingular point}. Our goal, as presented in Section~\ref{iden_singular}, is to locate the singular boundary points of nonlinearity intervals in $\interior(\mathcal{E})$, and then identify the transition points among the singular points, see Section~\ref{iden_transition}. 

\subsubsection{Computation of singular boundary points}\label{iden_singular}
Singular points of parameterized systems are well-studied in algebraic geometry, e.g., Sylvester's $19^{\rm th}$ century work in discriminants and resultants, see e.g.,~\cite{SJJ1851}. From a computational algebraic geometry viewpoint, the problem of computing singular boundary points for a parametric SDO problem was studied by the first and third authors in~\cite{HT} in a more general context. Here, we present a simplified process to 
locate the boundary points of nonlinearity intervals. Given an initial point $\bar{\epsilon} \in \interior(\mathcal{E})$ with a nonsingular Jacobian $J\!\big(V^*(\bar{\epsilon}),\bar{\epsilon}\big)$, the key idea is using Davidenko's~\cite{D1953,K1977} ordinary differential equation (ODE)
\begin{align}\label{eq:Davidenko}
J(V,\epsilon)\frac{dV}{d\epsilon}+\frac{\partial F(V,\epsilon)}{\partial \epsilon}=0
\end{align}
to track an optimal solution $V(\epsilon)$ from $\bar{\epsilon}$
to a boundary point in each direction. Since solutions of \eqref{eq:Davidenko} correspond to level sets of $F(V,\epsilon)$, i.e., $\{(V,\epsilon) : F(V,\epsilon)=c\}$ for arbitrary constant $c$, using the initial condition $V(\bar{\epsilon}) = V^*(\bar{\epsilon})$ yields the set of solutions to~\eqref{KKT_conditions} and~\eqref{eq:KKT_INEQ} for all $\epsilon$ in a neighborhood of $\bar{\epsilon}$. Hence, this approach utilizes the local information provided by the Jacobian, when it is nonsingular, to obtain accurate approximations of the optimal solutions nearby. The following theorem provides a summary of the solution~\cite{HT}.
\begin{theorem}\label{thm:V_ODE}
Let $\mathcal{I}_{\mathrm{reg}} \subseteq \interior(\mathcal{E})$ be an open interval containing $\bar{\epsilon}$ such that $J\!\big(V^*(\epsilon),\epsilon\big)$ is nonsingular for every $\epsilon \in \mathcal{I}_{\mathrm{reg}}$. Then, $V^*(\epsilon)$ is analytic on $\mathcal{I}_{\mathrm{reg}}$, and it is the unique solution~of
\begin{align}\label{ODE_system}
\frac{dV}{d\epsilon}=-J(V,\epsilon)^{-1}{\frac{\partial F(V,\epsilon)}{\partial \epsilon}}, \qquad V(\bar{\epsilon})= V^*(\bar{\epsilon}), \quad \epsilon \in \mathcal{I}_{\mathrm{reg}}.
\end{align}
\end{theorem}
\proof{Proof.}{See the Appendix.}
\endproof
\noindent
Using Theorem~\ref{thm:V_ODE} and the results of~\cite{HHL14}, we can track along $\mathcal{I}_{\mathrm{reg}}$, on which the optimal solution $V^*(\epsilon)$ is analytic by the implicit function theorem~\cite[Theorem~10.2.4]{D60}, until we reach the boundary points of $\mathcal{I}_{\mathrm{reg}}$. Thus, as the perturbation parameter approaches 
a singular boundary point of $\mathcal{I}_{\mathrm{reg}}$, ill-conditioning of $F(V,\epsilon)=0$, or spurious numerical behavior will be detected numerically. Consequently, we can avoid jumping over a transition point by using any reasonable mesh size that is sufficiently small for solving the ODE system in Theorem~\ref{thm:V_ODE}.
\begin{remark}
Theorem~\ref{thm:V_ODE} and the ODE system~\eqref{eq:Davidenko} serve as the basis of Algorithm~\ref{singular_identification} in Section~\ref{sec:partition_alg}.  \halmos
\end{remark}
\subsubsection{Identification of transition points}\label{iden_transition}
At a singular boundary point $\hat{\epsilon}$, we examine the uniqueness of the corresponding optimal solution $V^{a}(\hat{\epsilon})$, where $V^{a}(\hat{\epsilon})$ is an accumulation point of the sequence of unique optimal solutions $V^*(\epsilon)$, obtained from~\eqref{eq:Davidenko}, as $\epsilon \nearrow \hat{\epsilon}$ or $\epsilon \searrow \hat{\epsilon}$. An accumulation point exists, by the outer semicontinuity of $\mathcal{P}^*$ and $\mathcal{D}^*$ relative to $\interior(\mathcal{E})$, and it belongs to $\mathcal{P}^*(\hat{\epsilon}) \times \mathcal{D}^*(\hat{\epsilon})$. Toward this end, we compute the local dimension of the algebraic set $\mathbf{V}\big(F(V,\hat{\epsilon})\big)$ at $V^{a}(\hat{\epsilon})$ using a numerical local dimension test~\cite{BHPS09,WHS11}. The local dimension is defined as the maximum dimension of the irreducible components of $\mathbf{V}\big(F(V,\hat{\epsilon})\big)$, i.e., minimal algebraic subsets of $\mathbf{V}\big(F(V,\hat{\epsilon})\big)$, which contain $V^{a}(\hat{\epsilon})$, see Example~\ref{ex:local_dim}. A detailed description of algebraic sets and irreducible components can be found in~\cite{SW05}.

\vspace{5px}
\noindent
If $\mathbf{V}\big(F(V,\hat{\epsilon})\big)$ has 
local dimension zero at $V^{a}(\hat{\epsilon})$, then we can conclude from Lemma~\ref{continuity_nonlinearity} that $\hat{\epsilon}$ is a transition point,
since $V^a(\hat{\epsilon})$ turns out to be the unique optimal solution of $(\mathrm{P_{\hat{\epsilon}}})$ and $(\mathrm{D_{\hat{\epsilon}}})$. 
Otherwise, we need to examine the change of rank at a maximally complementary optimal solution $V^*(\hat{\epsilon})$. Such a solution is generic on the irreducible component of $\mathbf{V}\big(F(V,\hat{\epsilon})\big)$ which contains $V^{a}(\hat{\epsilon})$, and it can be computed efficiently using numerical algebraic geometry~\cite{BSHW2013}. 

\begin{example}\label{ex:local_dim}
For the system
\begin{align*}
F((x_1,x_2),\epsilon)=\begin{pmatrix} x_1^2+x_2^2-\epsilon\\(x_1^2+x_2^2-1)x_1 \end{pmatrix}
\end{align*}
the Jacobian with respect to $(x_1,x_2)$ is only singular at $\epsilon=0,1$. It is easy to see that $\mathbf{V}\big(F((x_1,x_2),0)\big)=\{(0,0)\}$ with a local dimension zero, while $\mathbf{V}\big(F((x_1,x_2),1)\big)=\{(x_1,x_2): x_1^2+x_2^2-1=0\}$ has local dimension one.  \halmos
\end{example}

\begin{remark}
The local dimension test serves as the basis of Algorithm~\ref{transition_identification} in Section~\ref{sec:partition_alg}. \halmos
\end{remark}
\subsubsection{Topology of singular points}
Although the set of transition points is always finite, in practice, the singular points need not be isolated. A case with infinitely many real singular points is demonstrated in Section~\ref{Exp:circle_line}, where every $V^*(\epsilon)$ in the only nonlinearity interval has a nonsingular Jacobian, see also Example~\ref{ex:singularity}. However, under the existence of a \textcolor{black}{generic} nonsingular point in $\interior(\mathcal{E})$, the algebraic formulation~\eqref{KKT_conditions} shows that the set of singular points must be an algebraic subset of $\mathbb{C}$, leading to the following finiteness result. 
\begin{proposition}\label{extrema2singular}
Assume that there exists a \textcolor{black}{generic} nonsingular point $\bar{\epsilon} \in \interior(\mathcal{E})$. 
Then the set of singular points in $\interior(\mathcal{E})$ is finite. As a consequence, the set of points at which $\mathcal{P}^*$ or $\mathcal{D}^*$ fails to be continuous relative to $\interior(\mathcal{E})$ is finite.
\end{proposition}
\proof{Proof.}
By definition, the set $\Upsilon$ of all $(\underline{V}(\epsilon),\epsilon)$ with a singular Jacobian satisfies
\begin{align}\label{singular_Vepsilon}
\Upsilon\!:=\!\big\{(V,\epsilon) \in \mathbb{C}^{m+2t(n)+1} : F(V,\epsilon) = 0, \ \det\!\big(J(V,\epsilon)\big) = 0 \big\}, 
\end{align}
where~\eqref{singular_Vepsilon} is a basic constructible set~\cite{BPR06} in $\mathbb{C}^{m+2t(n)+1}$. Since the projection of a constructible set to $\mathbb{C}$ is a constructible subset of $\mathbb{C}$~\cite[Theorem 1.22]{BPR06}, it holds that 
\begin{align}\label{singular_epsilon}
\big\{\epsilon \in \mathbb{C} : \exists V \in \mathbb{C}^{m+2t(n)} \ \st \ (V,\epsilon) \in \Upsilon \big\}
\end{align}
is either finite or the complement of a finite subset of $\mathbb{C}$, see e.g.,~\cite[Exercise 1.2]{BPR06}. On the other hand, it follows from the assumption and the implicit function theorem that the complement of~\eqref{singular_epsilon} contains an open neighborhood of $\bar{\epsilon}$. All this implies that the projection of $\Upsilon$ is finite, and thus it is an algebraic subset of $\mathbb{C}$. The finiteness result naturally holds when we restrict the set of singular points to $\mathbb{R}$, in which our domain $\mathcal{E}$ is defined. Consequently, there are only finitely many real singular points in $\interior(\mathcal{E})$.  \Halmos
\endproof

\begin{remark}
As a consequence of Proposition~\ref{extrema2singular} and~\cite[Theorem~5.12]{L13}, the polynomial system~\eqref{KKT_conditions} is zero-dimensional at every nonsingular $\epsilon \in \interior(\mathcal{E})$, i.e., $\mathbf{V}\big(F(V,\epsilon)\big)$ has only finitely many solutions almost everywhere on $\interior(\mathcal{E})$. \halmos
\end{remark}

\vspace{5px}
\noindent
The condition of Proposition~\ref{extrema2singular} is a global condition which requires that every solution of the algebraic set $\mathbf{V}\big(F(V,\epsilon)\big)$ at a \textcolor{black}{generic} $\epsilon \in \interior(\mathcal{E})$ has a nonsingular Jacobian. \textcolor{black}{Notice that $\mathbf{V}\big(F(V,\epsilon)\big)$ has a generic behavior over all $\epsilon \in \mathbb{C}$. In particular, there are only finitely many points $\mathcal{F} \subset \mathbb{C}$ which can have a different irreducible decomposition than the generic case. Hence, for any open interval $\mathcal{O} \subset \mathbb{R}$, there are at most finitely many points which are not generic. Therefore, $\epsilon \in \mathcal{O}$ is a generic nonsingular point if $\epsilon \not \in \mathcal{F}$ and every solution of $\mathbf{V}\big(F(V,\epsilon)\big)$ is nonsingular.}  

\vspace{5px}
\noindent
Recall from Lemma~\ref{nonsingularity_conditions} that strict complementarity and non-degeneracy conditions at $\epsilon$ are necessary and sufficient for the existence of a unique $V^*(\epsilon)$ with a nonsingular Jacobian. Therefore, the condition of Proposition~\ref{extrema2singular} is at least as strong as strict complementarity and non-degeneracy conditions. Interestingly, the following proposition indicates that for the polynomial system~\eqref{KKT_conditions} with generic data, there exists a nonsingular point $\bar{\epsilon}$ with probability~1.
\begin{proposition}
\textcolor{black}{The condition of Proposition~\ref{extrema2singular}} is a generic property with respect to all $(\mathcal{A}, b, C, \bar{C})$.
\end{proposition}
\proof{Proof.}
Without loss of generality, we will simply consider when $\epsilon=0$. It follows from~\cite[Theorem~7]{NRS10} that for generic $(\mathcal{A},b,C)$, all complex solutions of $F(V,0)=0$ are isolated and have nonsingular Jacobian. All this implies that for generic $(\mathcal{A},b,C)$, $\epsilon = 0$ is a nonsingular point. \halmos
\endproof

\begin{example}\label{ex:singularity}
There are special cases where the solution set $\mathbf{V}\big(F(V,\epsilon)\big)$
consists of isolated solutions or algebraic subsets with positive dimension. For instance, for the system
\begin{align*}
F((x_1,x_2),\epsilon) = \begin{pmatrix} (x_1^2+x_2^2-1)(x_1-x_2)\\ (x_1^2+x_2^2-1)(x_1-\epsilon)\end{pmatrix}
\end{align*}
there are two solution sets at $\epsilon\neq\pm\frac{1}{\sqrt{2}}$: a circle $\{(x_1,x_2) \in \mathbb{C}^2 : x_1^2+x_2^2=1\}$ and an isolated solution $(\epsilon,\epsilon)$. \halmos
\end{example}

\subsection{Partitioning algorithm}\label{sec:partition_alg}
Based on the descriptions in Sections~\ref{iden_singular} and~\ref{iden_transition} and the auxiliary problems in~\eqref{auxiliary_problems}, we present the outline of our numerical procedure, Algorithm~\ref{outermost_alg}. Algorithm~\ref{outermost_alg} consecutively calls the subroutines in Algorithms~\ref{invariancy_identification},~\ref{singular_identification}, and~\ref{transition_identification} to compute invariancy intervals, nonlinearity intervals, and transition points in $\interior(\mathcal{E})$. For the ease of exposition, see Remark~\ref{weaker_assumption}, we outline the pseudo codes by assuming, only in this section, the condition of Proposition~\ref{extrema2singular}. This condition will enable us to decompose $\interior(\mathcal{E})$ into the union of finitely many open intervals of maximal length by locating their finitely many singular boundary points.

\vspace{5px}
\noindent
In our numerical procedure, Algorithm~\ref{invariancy_identification} computes the boundary points of an invariancy interval by solving auxiliary problems~\eqref{auxiliary_problems} and then updates the set of transition points and the collection of invariancy intervals in $\interior(\mathcal{E})$. When Algorithm~\ref{invariancy_identification} fails to identify an invariancy interval, Algorithms~\ref{singular_identification} and~\ref{transition_identification} are subsequently called to locate the boundary points of a nonlinearity interval, if they exist, or to conclude the existence of a transition point. More specifically, this is done by locating the singular points in the remaining subinterval of $\interior(\mathcal{E})$, as described in Sections~\ref{iden_singular} and~\ref{iden_transition}:
\begin{itemize}
\item Algorithm~\ref{singular_identification} tracks the optimal solution of $(\mathrm{P_{\epsilon}})$ and $(\mathrm{D_{\epsilon}})$ by solving the ODE system~\eqref{eq:Davidenko} using a predictor-corrector tracking method~\cite{Butcher2003} until it detects a singular boundary point. 
 \item Algorithm~\ref{transition_identification} classifies singular points into transition and non-transition points. 
 \end{itemize}
\vspace{5px} 
Algorithm~\ref{singular_identification} is repeatedly called alongside Algorithm~\ref{invariancy_identification} until all invariancy intervals and singular points in $\interior(\mathcal{E})$ are identified. Finally, the collection of nonlinearity intervals are formed by removing the invariancy intervals and transition points from $\interior(\mathcal{E})$.

\vspace{5px}
\noindent
In order to completely cover the interval, the increment change $\Delta\epsilon$ can be positive or negative to allow both left and right movements from the starting point. Furthermore, we assume, for the simplicity of computation, that the domain $\mathcal{E}$ is bounded, i.e., $\mathcal{E}=[\mathcal{E}_{\min},\mathcal{E}_{\max}]$, where \hbox{$|\mathcal{E}_{\min}|,|\mathcal{E}_{\max}| < \infty$}. Accordingly, the optimal value of the auxiliary problems~\eqref{auxiliary_problems} is constrained to $[\mathcal{E}_{\min},\mathcal{E}_{\max}]$. For the sake of brevity, Algorithms~\ref{outermost_alg} through~\ref{transition_identification} only present the computation of invariancy intervals, nonlinearity intervals, and transition points on the subinterval $[\bar{\epsilon},\mathcal{E}_{\max})$, where $\bar{\epsilon}$ is the initial point.

\begin{remark}
Our approach is in direct contrast with finding transition points through solving $(\mathrm{P_{\epsilon}})-(\mathrm{D_{\epsilon}})$ on an arbitrarily meshed interval. In the latter case, as mentioned at the beginning of Section~\ref{nonlinearity_interval_algorithm}, only very refined mesh sizes may prevent the miscount of the transition points. \halmos
\end{remark}

\paragraph{Computation of singular points and invariancy intervals}
Theorem~\ref{thm:V_ODE} specifies a systematic way to approximate the boundary points of the interval $\mathcal{I}_{\mathrm{reg}}$ surrounding the given~$\bar{\epsilon}$. The numerical detection of singular points is described in detail in~\cite{HT} with respect to several singularity criteria, e.g., the derivative of $\lambda_{\min}\big(X^*(\epsilon)\big)$ and $\lambda_{\min}\big(S^*(\epsilon)\big)$ with respect to $\epsilon$, or the singularity of the Jacobian of~\eqref{KKT_conditions}. We omit the details here and refer the reader to~\cite{HT} for more information on the numerical implementation of the singularity criteria. 

\vspace{5px}
\noindent
Once a singular point is identified, the numerical solution obtained from the ODE system~\eqref{eq:Davidenko} at the next mesh point is most likely non-optimal, due to the numerical instability or the infeasibility of the solution. Thus, we invoke a primal-dual IPM in Algorithms~\ref{invariancy_identification} and~\ref{singular_identification} to compute the unique optimal solution at the first neighboring mesh point in the remaining interval. In order to guarantee that every singular point is correctly identified, a finer mesh pattern might be needed, and a higher precision might be required for the computation of singular points, far beyond the double precision arithmetic.

\paragraph{Solution sharpening}
The process of increasing the algebraic precision of a singular point is known as the sharpening process, see Algorithm~\ref{singular_identification}. Since the singular points are algebraic numbers, they can be computed to arbitrary accuracy, see e.g.,~\cite{HS17}. 
 More specifically, using a numerical approximation of a given singular point, which is indeed the nearest mesh point to the singular point, the theory of isosingular sets~\cite{HW13} allows one to construct a new polynomial system where Newton's method would converge quadratically to the singular point.  
 
\paragraph{Classification of singular points}
 The use of adaptive precision, see e.g.,~\cite{BHSW08},
 in {\tt Bertini}~\cite{BHSW06,BSHW2013} ensures that adequate precision is being used
 for reliable computations near the singular solutions. This method enables one to compute a maximally complementary optimal solution near $V^{a}(\hat{\epsilon})$ to arbitrary accuracy. With the ability to refine the accuracy of a maximally complementary optimal solution, we can determine if a given singular point is a transition point. This can be done robustly by examining the rank of $X^*(\epsilon)$ and $S^*(\epsilon)$ using standard numerical rank revealing methods, such as singular value decomposition. More specifically, by computing the eigenvalues of an approximate maximally complementary optimal solution at various precisions, one can determine if the least positive eigenvalues of $X^*(\epsilon)$ and $S^*(\epsilon)$ converge to zero as we increase the precision of computation. This process accurately reveals the rank of $X^*(\epsilon)$ and $S^*(\epsilon)$ at a singular point.

\begin{remark}\label{weaker_assumption}
The sole purpose of imposing the condition of Proposition~\ref{extrema2singular} in Algorithm~\ref{outermost_alg} is to ensure finite decomposition of $\interior(\mathcal{E})$. Otherwise, Algorithm~\ref{singular_identification} can be individually applied to find a subinterval of the nonlinearity interval, even under a weaker condition than Proposition~\ref{extrema2singular}. More precisely, the existence of $\bar{\epsilon}$ with a nonsingular $J\big(V^*(\bar{\epsilon}),\bar{\epsilon}\big)$ is all we need in Theorem~\ref{thm:V_ODE} to compute a subinterval of a nonlinearity interval containing $\bar{\epsilon}$, see the proof of Theorem~\ref{thm:V_ODE} in the Appendix. Without the condition of Proposition~\ref{extrema2singular}, however, a full decomposition of $\interior(\mathcal{E})$ may not be possible using Algorithm~\ref{outermost_alg}, because singular points need not be isolated in that case. \halmos
\end{remark}

\begin{algorithm}[]
\small
\caption{Partitioning of $\interior(\mathcal{E})$}
\label{outermost_alg}
\begin{algorithmic}[l]
\State\textbf{Global Input:} Problem data: $\mathcal{A}$, $b$, $C$, $\bar{C}$, and the domain  $\mathcal{E}=[\mathcal{E}_{\min},\mathcal{E}_{\max}]$.
\State\textbf{Local Input:} an initial point $\epsilon_{\mathrm{init}}\in\interior(\mathcal{E})$ with a nonsingular Jacobian $J\!\big(V^*(\epsilon_{\mathrm{init}}),\epsilon_{\mathrm{init}}\big)$, a positive increment change $\Delta\epsilon$.
\State \textbf{Output:} $\mathcal{U}_{\mathrm{inv}}$: union of invariancy intervals in $(\mathcal{E}_{\min},\mathcal{E}_{\max})$,\\
\qquad\qquad \ $\mathcal{U}_{\mathrm{non}}$:  union of nonlinearity intervals in $(\mathcal{E}_{\min},\mathcal{E}_{\max})$,\\ 
\qquad\qquad \  $\mathcal{U}_{\mathrm{tran}}$: set of transition points in $(\mathcal{E}_{\min},\mathcal{E}_{\max})$.

\vspace{15px}
\noindent
\textbf{Procedure:}
\vspace{10px}
\begin{itemize}
\item Set $\epsilon = \epsilon_{\mathrm{init}}$, $\mathcal{U}_{\mathrm{inv}}=\emptyset$, $\mathcal{U}_{\mathrm{non}}=(\mathcal{E}_{\min},\mathcal{E}_{\max})$, $\mathcal{U}_{\mathrm{tran}} = \emptyset$, and $\mathcal{U}_{\mathrm{sin}}=\emptyset$.
\end{itemize}

\vspace{10px}
\While{$\epsilon < \mathcal{E}_{\max}$} 	
\Repeat \Comment{Compute invariancy intervals}	
\State \textbullet \ Find invariancy intervals: Apply Algorithm~\ref{invariancy_identification} using the input $\Delta \epsilon$, $\epsilon$, $\mathcal{U}_{\mathrm{inv}}$, $\mathcal{U}_{\mathrm{non}}$, and $\mathcal{U}_{\mathrm{tran}}$
\State \ \ (Algorithm \ref{invariancy_identification} outputs $\alpha_{\mathrm{inv}}$ and $\beta_{\mathrm{inv}}$ and updates input arguments $\epsilon$, $\mathcal{U}_{\mathrm{inv}}$, $\mathcal{U}_{\mathrm{non}}$, and $\mathcal{U}_{\mathrm{tran}}$).
\Until{$\alpha_{\mathrm{inv}} < \beta_{\mathrm{inv}}$ and $\epsilon < \mathcal{E}_{\max}$} \\
		
\If{$\epsilon <\mathcal{E}_{\max}$} 	\Comment{Compute singular points}	
\State \textbullet \ Apply Algorithm~\ref{singular_identification} using the input $\Delta \epsilon$, $\epsilon$, $\mathcal{U}_{\mathrm{inv}}$, $\mathcal{U}_{\mathrm{non}}$, $\mathcal{U}_{\mathrm{sin}}$, and $\mathcal{U}_{\mathrm{tran}}$.                
\EndIf
\EndWhile	
\vspace{5px}	
\begin{itemize}	
\item Apply Algorithm~\ref{transition_identification} using the input $\mathcal{U}_{\mathrm{sin}}$ and $\mathcal{U}_{\mathrm{tran}}$.
\item Set $\mathcal{U}_{\mathrm{non}}=\mathcal{U}_{\mathrm{non}} \setminus \mathcal{U}_{\mathrm{tran}}$.  \Comment{Form the nonlinearity intervals}
\end{itemize}
\end{algorithmic}
\end{algorithm}

\begin{algorithm}[]
\small
\caption{Computation of invariancy intervals}
\label{invariancy_identification}
\begin{algorithmic}[l]
\State \textbf{Global Input:} Problem data: $\mathcal{A}$, $b$, $C$, $\bar{C}$, and the domain $\mathcal{E}=[\mathcal{E}_{\min},\mathcal{E}_{\max}]$.\\
\textbf{Local Input:} an increment change $\Delta \epsilon$, $\epsilon$, $\mathcal{U}_{\mathrm{inv}}$, $\mathcal{U}_{\mathrm{non}}$, $\mathcal{U}_{\mathrm{tran}}$.
\State \textbf{Output:} $(\alpha_{\mathrm{inv}}, \beta_{\mathrm{inv}})$ and updated $\epsilon$, $\mathcal{U}_{\mathrm{inv}}$, $\mathcal{U}_{\mathrm{non}}$, $\mathcal{U}_{\mathrm{tran}}$.\\
		
\vspace{15px}
\noindent
\textbf{Procedure:}
\vspace{10px}
\begin{itemize}
\item  Compute the unique optimal solution $V^*(\epsilon)$ using a primal-dual IPM.
\item  Compute the orthonormal basis $Q_{\mathcal{N}(\epsilon)}$ from $V^*(\epsilon)$.
\item Using $Q_{\mathcal{N}(\epsilon)}$ solve the pair of SDO problems~\eqref{auxiliary_problems} restricted to $[\mathcal{E}_{\min}, \mathcal{E}_{\max}]$ to compute the boundary points $\alpha_{\mathrm{inv}}$ and $\beta_{\mathrm{inv}}$.
\end{itemize}

\vspace{10px}
\If{$\alpha_{\mathrm{inv}} < \epsilon < \beta_{\mathrm{inv}}$} \Comment{An invariancy interval exists}
\begin{itemize}
\item Update the union of invariancy intervals by adding the newly found interval $(\alpha_{\mathrm{inv}},\beta_{\mathrm{inv}})$ to the union of invariancy intervals $\mathcal{U}_{\mathrm{inv}}$: $\mathcal{U}_{\mathrm{inv}}=\mathcal{U}_{\mathrm{inv}} \cup (\alpha_{\mathrm{inv}},\beta_{\mathrm{inv}})$. 
\item Update the union of nonlinearity intervals by removing the invariancy interval $(\alpha_{\mathrm{inv}},\beta_{\mathrm{inv}})$ from the current union of nonlinearity intervals: $\mathcal{U}_{\mathrm{non}}=\mathcal{U}_{\mathrm{non}} \setminus (\alpha_{\mathrm{inv}},\beta_{\mathrm{inv}})$. 
\item Update the set of transition points by 
\begin{align*}
\mathcal{U}_{\mathrm{tran}}=\begin{cases} \mathcal{U}_{\mathrm{tran}} \cup \{\alpha_{\mathrm{inv}}\} \ \  \alpha_{\mathrm{inv}} > \mathcal{E}_{\min},\\ \mathcal{U}_{\mathrm{tran}} \cup \{\beta_{\mathrm{inv}}\} \ \  \beta_{\mathrm{inv}} < \mathcal{E}_{\max}. \end{cases}
\end{align*} 
\vspace{5px}
\item Move past a transition point by $\epsilon = \beta_{\mathrm{inv}} + \Delta \epsilon$.	 
\end{itemize}
\EndIf	
\end{algorithmic}
\end{algorithm}	

\begin{algorithm}[]
\small
\caption{Computation of the singular points}
\label{singular_identification}
\begin{algorithmic}[l]
\State \textbf{Global Input:} Problem data: $\mathcal{A}$, $b$, $C$, $\bar{C}$, and the domain $[\mathcal{E}_{\min},\mathcal{E}_{\max}]$.
\State \textbf{Local Input:} $\Delta\epsilon$, $\epsilon$, $\mathcal{U}_{\mathrm{inv}}$, $\mathcal{U}_{\mathrm{non}}$, $\mathcal{U}_{\mathrm{sin}}$, $\mathcal{U}_{\mathrm{tran}}$.
\State \textbf{Output:}  Updated $\epsilon$, $\mathcal{U}_{\mathrm{inv}}$, $\mathcal{U}_{\mathrm{non}}$, $\mathcal{U}_{\mathrm{sin}}$, and $\mathcal{U}_{\mathrm{tran}}$.\\

\vspace{10px}
\noindent
\textbf{Procedure:}
\begin{itemize}
\item Compute the unique optimal solution $V^*:=V^*(\epsilon)$ using a primal-dual IPM.
\end{itemize}

\vspace{15px}
\While{Jacobian is nonsingular on $[\epsilon,\epsilon+\Delta\epsilon]$ and $\epsilon +\Delta\epsilon \in (\mathcal{E}_{\min},\mathcal{E}_{\max})$} \Comment{Check the singularity}
\begin{itemize}
\item Proceed to the next mesh point by $\epsilon=\epsilon+\Delta\epsilon$.
\item Compute the unique optimal solution $V^*(\epsilon)$ by solving~\eqref{eq:Davidenko} with the initial point $V^*$.
\end{itemize}
\EndWhile\\
	
\If{a singular point exists in $[\epsilon,\epsilon+\Delta\epsilon]$ and $\epsilon +\Delta\epsilon \in (\mathcal{E}_{\min},\mathcal{E}_{\max})$} \Comment{A singular point exists}
\begin{itemize}
\item Use solution sharpening to compute the singular point $\hat{\epsilon}$ and set $\mathcal{U}_{\mathrm{sin}}=\mathcal{U}_{\mathrm{sin}}\cup \{(V^a(\hat{\epsilon}),\hat{\epsilon})\}$.
\item Move past the singular point by $\epsilon=\hat{\epsilon}+\Delta\epsilon$. 
\end{itemize}
\Else 
\begin{itemize}
\item Proceed to the next mesh point by $\epsilon=\epsilon+\Delta\epsilon$.
\end{itemize}
\EndIf
\end{algorithmic}
\end{algorithm}	

\begin{algorithm}[]
\small
\caption{Classification of the singular points}
\label{transition_identification}
\begin{algorithmic}[l]
 \State \textbf{Global Input:} Problem data: $\mathcal{A}$, $b$, $C$, and $\bar{C}$.
 \State \textbf{Local Input:} $\mathcal{U}_{\mathrm{sin}}$ and $\mathcal{U}_{\mathrm{tran}}$.
\State \textbf{Output:} Updated $\mathcal{U}_{\mathrm{tran}}$.
		
\vspace{15px}
\noindent
\textbf{Procedure:}
\vspace{10px}
\For{$(V,\epsilon) \in \mathcal{U}_{\mathrm{sin}}$}
\begin{itemize}
\item Calculate the local dimension $d$ of the algebraic set $\mathbf{V}\big(F(V,\epsilon)\big)$, defined in~\eqref{algebraic_set_optimal}, at $V$. 
\end{itemize}
\vspace{10px}
\If{$d=0$} \Comment{A transition point exists}
\State \textbullet \ Update the set of transition points by $\mathcal{U}_{\mathrm{tran}}=\mathcal{U}_{\mathrm{tran}} \cup \{\epsilon\}$.
\Else 	
		
\State \textbullet \ Use a polynomial solver to compute $V^*(\epsilon)$ in the irreducible component which contains $V$.
\vspace{10px}
\If{ the rank of $X^*(\epsilon)$ or $S^*(\epsilon)$ changes} \Comment{A transition point exists}
\State \textbullet \ Update the set of transition points by $\mathcal{U}_{\mathrm{tran}}=\mathcal{U}_{\mathrm{tran}} \cup \{\epsilon\}$.
\EndIf
\EndIf
\EndFor
\end{algorithmic}
\end{algorithm}
\section{Numerical examples}\label{numerical_experiments}
In this section, using the approaches described in Section~\ref{sec:partition_alg} and outlined by Algorithms~\ref{outermost_alg} through~\ref{transition_identification}, we conduct numerical experiments on the computation of invariancy intervals, nonlinearity intervals, and transition points. Section~\ref{Exp:circle_line} demonstrates the convergence rate of computing the singular boundary points. Section~\ref{Exp:Exp_ce} describes a parametric SDO problem where the continuity of the dual optimal set mapping fails at a transition point. Section~\ref{sec:non-transition} computes the nonlinearity interval of the parametric SDO problem~\eqref{3elliptope_cut_eq} where the Jacobian is singular at a non-transition point. All numerical experiments are conducted on a PC with Intel Core i7-6500U CPU @2.5 GHz.

\subsection{Convergence rate}\label{Exp:circle_line}
\begin{figure}
 \begin{minipage}[c]{0.3\textwidth}
\includegraphics[height=1.9in]{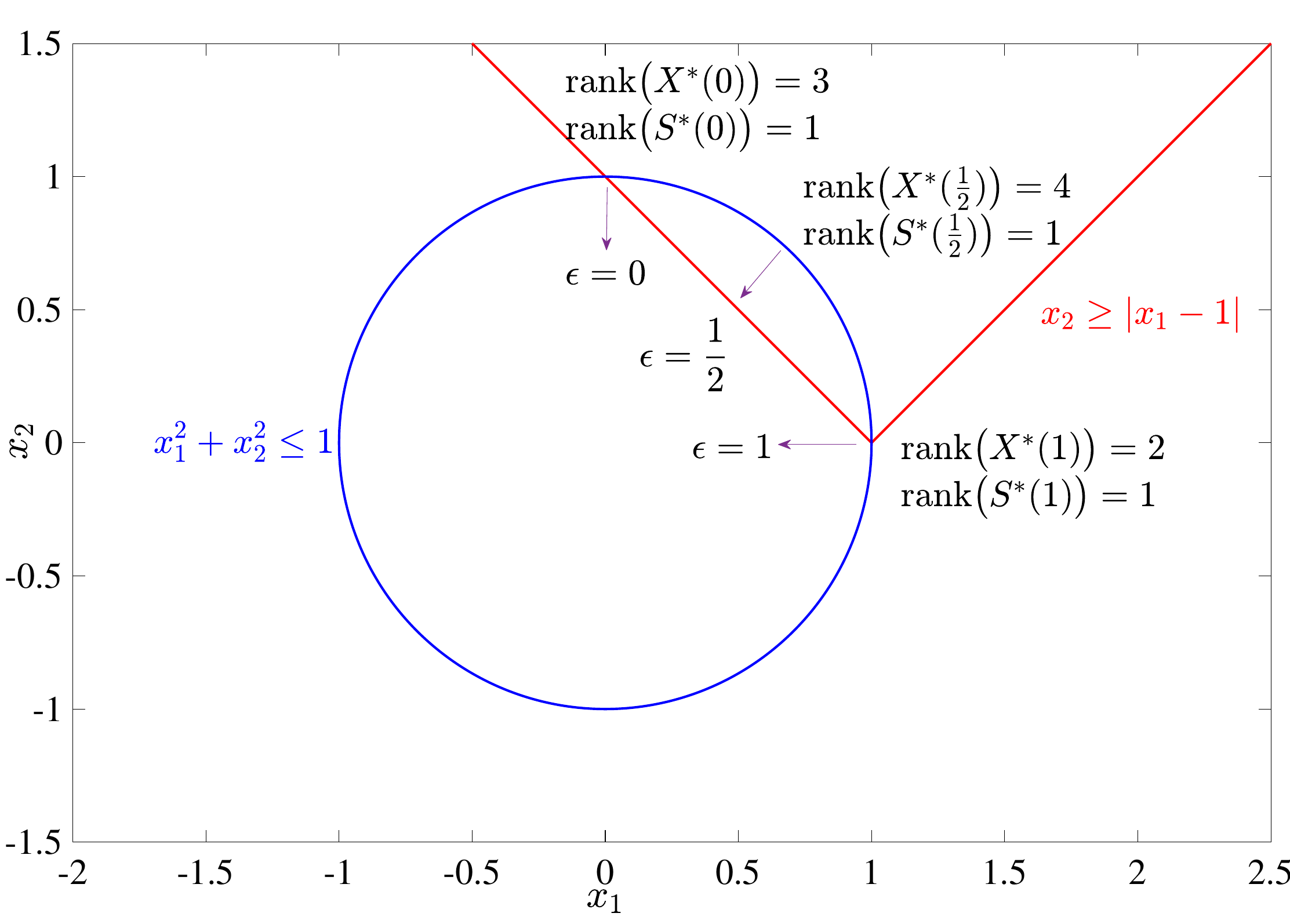}
\end{minipage}\hfill
\begin{minipage}[c]{0.7\textwidth}
\caption{The feasible set of problem \eqref{Eq:Exp_cl}.}
\label{Fig:Exp_ce}
\end{minipage}
\end{figure}

Consider the following parametric convex optimization problem
\begin{equation}\label{Eq:Exp_cl}
\begin{aligned}
\min \ \ &-2\epsilon x_1 - 2(1-\epsilon) x_2\\
\st \ \ & \begin{pmatrix}1 &x_1	& x_2 & 0 & 0\\x_1 & 1& 0 & 0& 0\\x_2 &0 &1& 0 & 0\\
0 &0 & 0 & x_2&x_1-1\\
0 & 0 & 0 &x_1-1&x_2\\
\end{pmatrix} \succeq 0,
\end{aligned}
\end{equation}
which can be cast into the primal form $(\mathrm{P_{\epsilon}})$, where $m=13$ and $X \in \mathbb{S}^5$. The block structure of the matrix indicates that~\eqref{Eq:Exp_cl} is indeed an SDO reformulation of a parametric second-order conic optimization problem with $\mathcal{E}=\mathbb{R}$, see also Figure~\ref{Fig:Exp_ce}. For computational purposes, we choose a bounded domain $[-\frac{1}{4},\frac{5}{4}]$ and the initial point $\epsilon = \frac14$, where $\rank\!\big(X^*(\frac14)\big)=4$, $\rank\!\big(S^*(\frac14)\big)=1$, and $J\big(V^*(\frac14),\frac14\big)$ is nonsingular. 

\begin{figure}
 \begin{minipage}[c]{0.6\textwidth}
\includegraphics[height=2.8in]{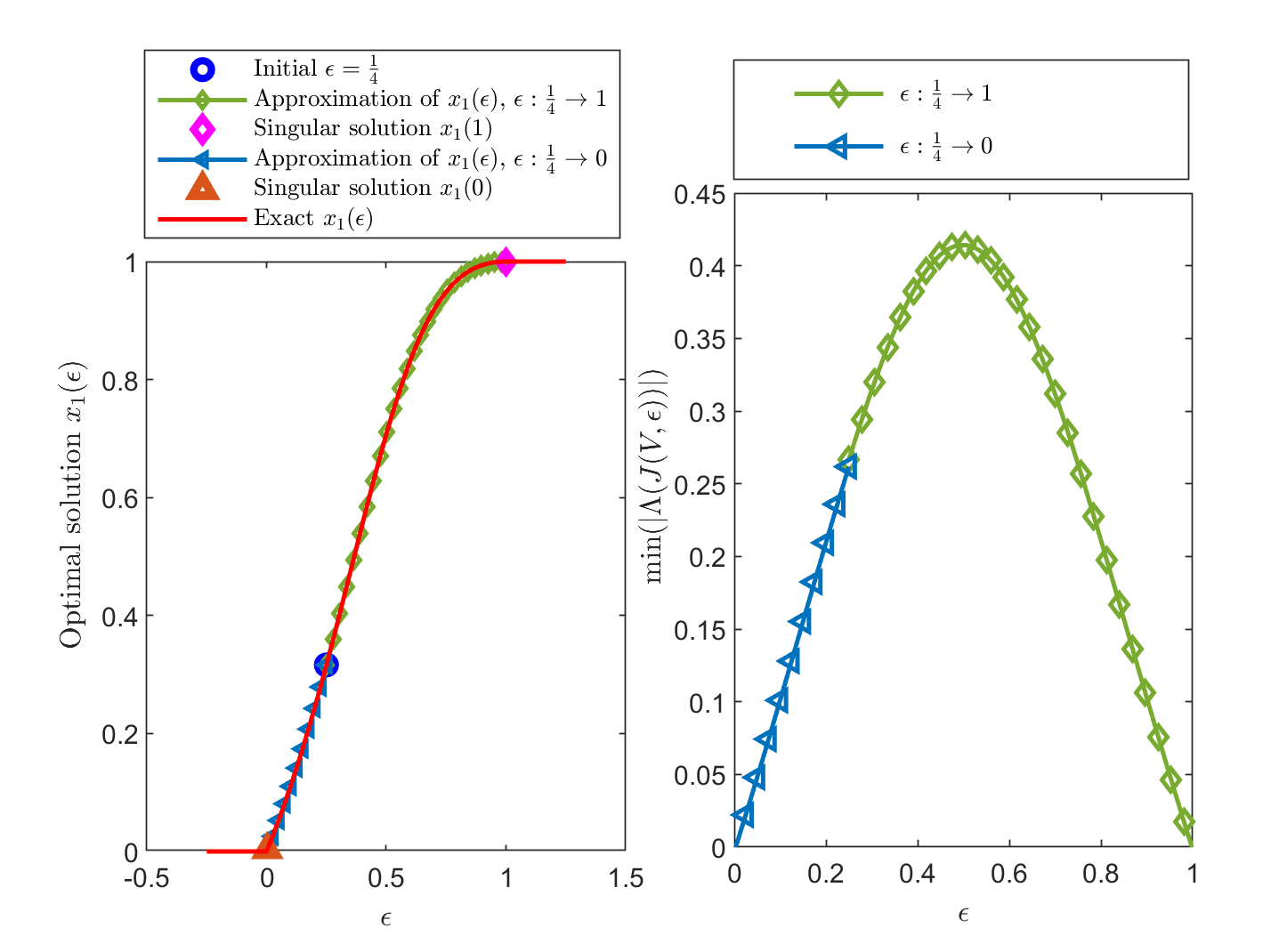}
\end{minipage}\hfill
\begin{minipage}[c]{0.4\textwidth}
\caption{Left: The exact and numerical approximation of $x_1(\epsilon)$ versus $\epsilon$. Right: The minimum modulus of the Jacobian eigenvalues.}
\label{Fig:Exp_cl_num}
\end{minipage}
\end{figure}

\vspace{5px}
\noindent 
Algorithm~\ref{invariancy_identification} identifies $\epsilon= \frac14$ as a point belonging to a nonlinearity interval. We then invoke Algorithm~\ref{singular_identification} to track the unique optimal solutions until we locate the boundary points $\epsilon=0$ and $\epsilon=1$. Algorithm~\ref{singular_identification} then computes a sufficiently accurate approximation of the boundary points. Figure~\ref{Fig:Exp_cl_num} demonstrates the exact and numerical approximation of
$x_1(\epsilon)$ and the minimum modulus of the Jacobian eigenvalues versus $\epsilon$. In particular, this tracking indicates that the Jacobian approaches singularity near $\epsilon=0$ and $\epsilon=1$.  

\vspace{5px}
\noindent
Restarting at the first mesh point next to the boundary points, Algorithm~\ref{invariancy_identification} identifies the invariancy intervals $(-\frac14,0)$ and $(1,\frac54)$ and determines that $\epsilon=0$ and $\epsilon=1$ are indeed the transition points of the optimal partition.

\vspace{5px}
\noindent
We point out that the condition of Proposition~\ref{extrema2singular} fails in this case. More specifically, for every $\epsilon \in \mathbb{R}$ the block diagonal structure in~\eqref{Eq:Exp_cl} allows for infinitely many real solutions $\underline{V}(\epsilon)=\big(\svectorize\!\big(\underline{X}(\epsilon)\big);\underline{y}(\epsilon);\svectorize\!\big(\underline{S}(\epsilon)\big)\big)$ for~\eqref{KKT_conditions}, such that 
\begin{align*}
 \underline{X}(\epsilon)=\begin{pmatrix} 1 & 1 & 0 & 0 & 0 \\1 & 1 & 0 & 0 & 0\\0 & 0 & 1 & 0 & 0\\0 & 0 & 0 & 0 & 0\\0 & 0 & 0 & 0 & 0 \end{pmatrix}, \quad   \underline{S}(\epsilon)=\begin{pmatrix} \ \ \epsilon + \zeta & -\epsilon -\zeta & 0 & 0 & 0 \\ -\epsilon - \zeta & \ \ \epsilon + \zeta & 0 & 0 & 0 \\ \ \ 0 & \ \ 0 & 0 & 0 & 0\\ \ \ 0 & \ \ 0 & 0 & 2(\epsilon-1) & \zeta\\ \ \ 0 & \ \ 0 & 0& \zeta & 0 \end{pmatrix}, \qquad \forall \zeta \in \mathbb{R}.
 \end{align*}
\noindent
Nevertheless, since the Jacobian $J\big(V^*(\frac14),\frac14\big)$ is nonsingular, the weaker condition described in 
Remark~\ref{weaker_assumption} holds and thus Algorithm~\ref{singular_identification} still correctly produces the boundary points of the nonlinearity interval.

\vspace{5px}
\noindent	
Using different patterns of mesh points, we demonstrate the convergence of $x_1(\epsilon)$, computed by Algorithm~\ref{singular_identification}, when $\epsilon$ approaches the singular boundary points $\epsilon=0$ and $\epsilon =1$. To that end, we let initial $\Delta\epsilon$ take values from $0.05\times 2^{-j}$ for $j = 0,\ldots, 5$ or $0.03\times 2^{-j}$ for $j = 0,\ldots, 5$, and we set $\epsilon = \frac14$ as the initial point. Tables~\ref{Tab:Exp3_convergence1} and~\ref{Tab:Exp3_convergence2} summarize the numerical results, where the $L_1$ error between the exact and numerical approximation of $x_1(\epsilon)$ on $[\frac14,1)$ and $(0,\frac14]$, the order of convergence, and the computation time are reported. The order of convergence is computed by
\begin{align*}
\rho_{j+1}\!:=\!\log_2\Bigg(\frac{\mathrm{Err}(\Delta\epsilon_{j})}{\mathrm{Err}(\Delta\epsilon_{j+1})}\Bigg), \qquad j=0,\ldots, 4,
\end{align*}
where $\mathrm{Err}(\Delta\epsilon_j)$ denotes the $L_1$ error associated with mesh pattern $j$. Notice the difference between $\rho_j$ and the classical notion of the order of convergence in computational optimization.
\begin{table}[H]
\centering
\small
\caption{Convergence of $x_1(\epsilon)$ when $\epsilon$ approaches the singular point $\epsilon = 1$.}
\label{Tab:Exp3_convergence1}
\begin{tabular}{c|c|c|c|c||cc}
\hline
$j$&\multicolumn{1}{c|}{$\Delta\epsilon_j$}&Approximate singular point & $\mathrm{Err}(\Delta\epsilon_j)$ & $\rho_j$&\multicolumn{1}{c}{CPU(s)}\\\hline 		
0 & $0.05$&1.00&$4.1597\times 10^{-6}$&-&4.05\\
1 & $0.05\times 2^{-1}$&1.00&$2.6520\times 10^{-7}$&3.971&6.56\\ [0.03in]
2 & $0.05\times 2^{-2}$&1.00&$1.6707\times 10^{-8}$&3.989&12.79\\ [0.03in]
3 & $0.05\times 2^{-3}$&1.00&$1.0484\times 10^{-9}$&3.994&26.14\\ [0.03in]
4 & $0.05\times 2^{-4}$&1.00&$6.5671\times 10^{-11}$&3.997&55.81\\ [0.03in]
5 & $0.05\times 2^{-5}$&1.00&$4.1090\times 10^{-12}$&3.998&125.27\\ [0.03in]
\hline
\end{tabular}
\end{table}

\vspace{5px}
\noindent	
In Table~\ref{Tab:Exp3_convergence1}, the singular point $\epsilon = 1$ is exactly identified by Algorithm~\ref{singular_identification}, since the singular point coincides with one of the mesh points. In general, however, it is unlikely that a singular point belongs to the mesh point set. This can be observed in Table~\ref{Tab:Exp3_convergence2}, where a fixed increment change $0.03\times 2^{-j}$ for $j = 0,\ldots, 5$ is utilized.  
In this case, the approximate singular point is taken as the
last mesh point before the minimum eigenvalues of $X^*(\epsilon)$ or $S^*(\epsilon)$, obtained from the ODE system~\eqref{eq:Davidenko}, become negative, or the first mesh point at which the minimum modulus of the Jacobian eigenvalues drops below $10^{-5}$.  
As stated in Section~\ref{sec:partition_alg}, we can utilize numerical algebraic geometric
tools to compute a singular point to arbitrary accuracy, but at the expense of increasing computational time.

\begin{table}[H]
\centering
\small
\caption{Convergence of $x_1(\epsilon)$ when $\epsilon$ approaches the singular point $\epsilon = 0$.}
\label{Tab:Exp3_convergence2}
\begin{tabular}{c|c|c|c|c||cc}
\hline
$j$&\multicolumn{1}{c|}{$\Delta\epsilon_j$}& \multicolumn{1}{c|}{Approximate singular point} & $\mathrm{Err}(\Delta\epsilon_j)$ & $\rho_j$&\multicolumn{1}{c}{CPU(s)}\\\hline 	0 & $0.03$&0.01&$2.0415\times 10^{-7}$&-& 2.85\\
1 & $0.03\times 2^{-1}$&0.01&$1.2917\times10^{-8}$&3.982&4.57\\ [0.03in]
2 & $0.03\times 2^{-2}$&0.025&$8.2444\times10^{-10}$&3.970&8.52\\ [0.03in]
3 & $0.03\times 2^{-3}$&0.0025&$5.1677\times10^{-11}$&3.996&17.73\\ [0.03in]
4 & $0.03\times 2^{-4}$&$0.000625$&$3.2461\times10^{-12}$&3.993&34.90\\ [0.03in]
5 & $0.03\times 2^{-5}$&$0.000625$&$2.0302\times10^{-13}$&3.999&72.34\\ [0.03in]
\hline
\end{tabular}
\end{table}

\subsection{A transition point with discontinuous dual optimal set mapping}\label{Exp:Exp_ce}
We next consider the parametric convex optimization problem
\begin{equation}\label{Eq:ellipse_circle}
\begin{aligned}
\min \ \ &\epsilon x_1 + (1-\epsilon) x_2\\
\st \ \  & \begin{pmatrix} 1 &x_1	& x_2 &0 &0 &0\\x_1 & 1& 0& 0& 0& 0\\x_2 &0 &1& 0& 0& 0\\
	0&0&0&1&\frac{1}{2}x_1&x_2\\
	0&0&0&\frac{1}{2}x_1&1&0\\
	0&0&0&x_2&0&1
\end{pmatrix} \succeq 0,
\end{aligned}
\end{equation}
in which the feasible set is compact and $\mathcal{E}=\mathbb{R}$. Analogous to~\eqref{Eq:Exp_cl}, this parametric problem can be cast into the primal form $(\mathrm{P_{\epsilon}})$ with $m=19$ and $X \in \mathbb{S}^6$. It can be verified that $J\!\big(V^*(\epsilon),\epsilon\big)$ is nonsingular, $\rank\!\big(X^*(\epsilon)\big)=5$, and $\rank\!\big(S^*(\epsilon)\big)=1$ at every $\epsilon \in \mathcal{E}\setminus\{0\}$. Since both the primal and dual problems have unique optimal solutions for every $\epsilon \in \mathcal{E}\setminus\{0\}$, the dual optimal set mapping fails to be continuous at $\epsilon = 0$.

\begin{figure}
\centering
\includegraphics[height=1.9in]{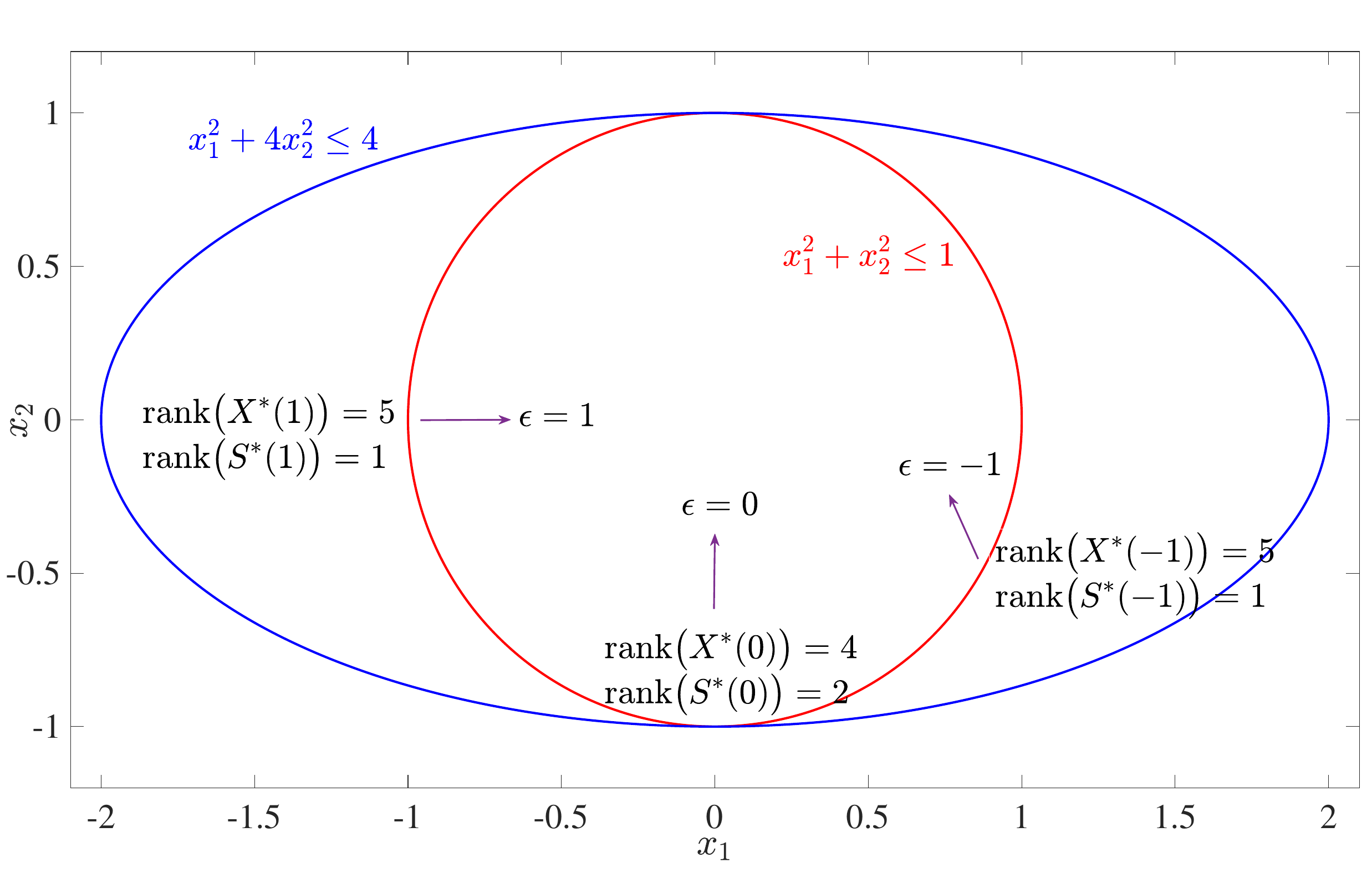}
\includegraphics[height=2.5in]{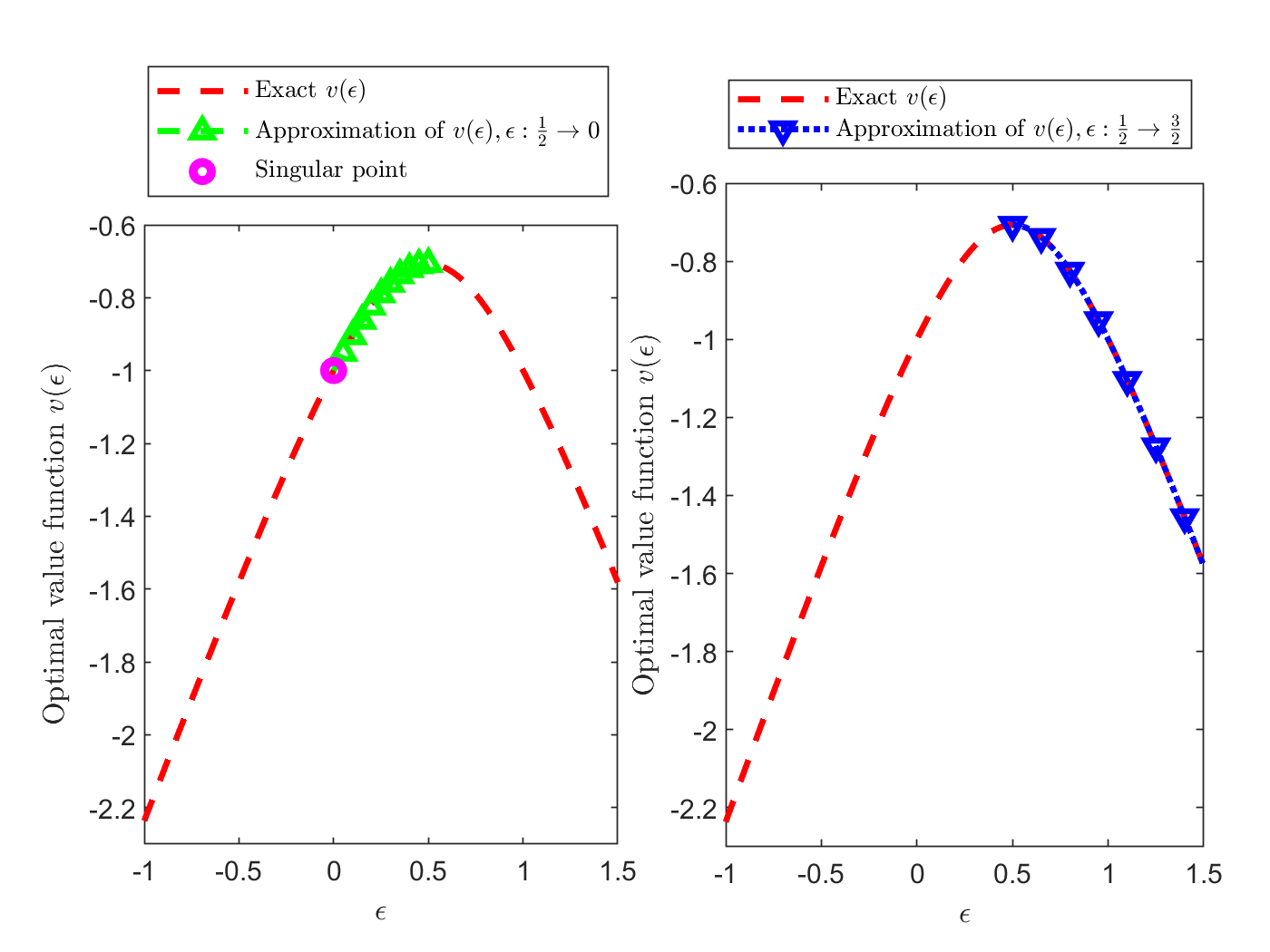}
\caption{Left: The feasible set of problem~\eqref{Eq:ellipse_circle}. Right: The exact and numerical approximation of the optimal value function for problem~\eqref{Eq:ellipse_circle} on $[-1,\frac32]$.}
\label{Fig:interval_ce}
\end{figure}

\vspace{5px}
\noindent
For the purpose of numerical experiments, we consider the bounded domain $[-1,\frac{3}{2}]$. When starting from initial point $\epsilon = \frac12$ with a fixed increment change $0.01$, Algorithm~\ref{singular_identification} properly identifies $\epsilon = 0$ as a singular boundary point. Figure~\ref{Fig:interval_ce} demonstrates the exact optimal value function versus its numerical approximation obtained from Algorithm~\ref{singular_identification}. 
Upon refining the accuracy of the approximate singular point and obtaining the singular point $\epsilon = 0$, we invoke {\tt Bertini} solver in Algorithm~\ref{transition_identification} to compute the dimension of all irreducible components of $\mathbf{V}\big(F(V,0)\big)$ which contain $V^a(0)$.
We observe that $V^a(0)$ lies on a 1-dimensional irreducible component of $\mathbf{V}\big(F(V,0)\big)$, and there exists a generic solution $V^*(0)$ such that $\rank\!\big(X^*(0)\big)=4$ and $\rank\!\big(S^*(0)\big)=2$. All this indicates that the rank of $X^*(\epsilon)$ and $S^*(\epsilon)$ change at $\epsilon = 0$, and thus $\epsilon = 0$ is a transition point. Consequently, we can partition $(-1,\frac{3}{2})$ into two nonlinearity intervals $(-1,0)$ and $(0,\frac{3}{2})$ and the transition point $\{0\}$.
\subsection{A non-transition point with singular Jacobian}\label{sec:non-transition}
\begin{figure}
 \begin{minipage}[c]{0.6\textwidth}
\includegraphics[width=3.5in]{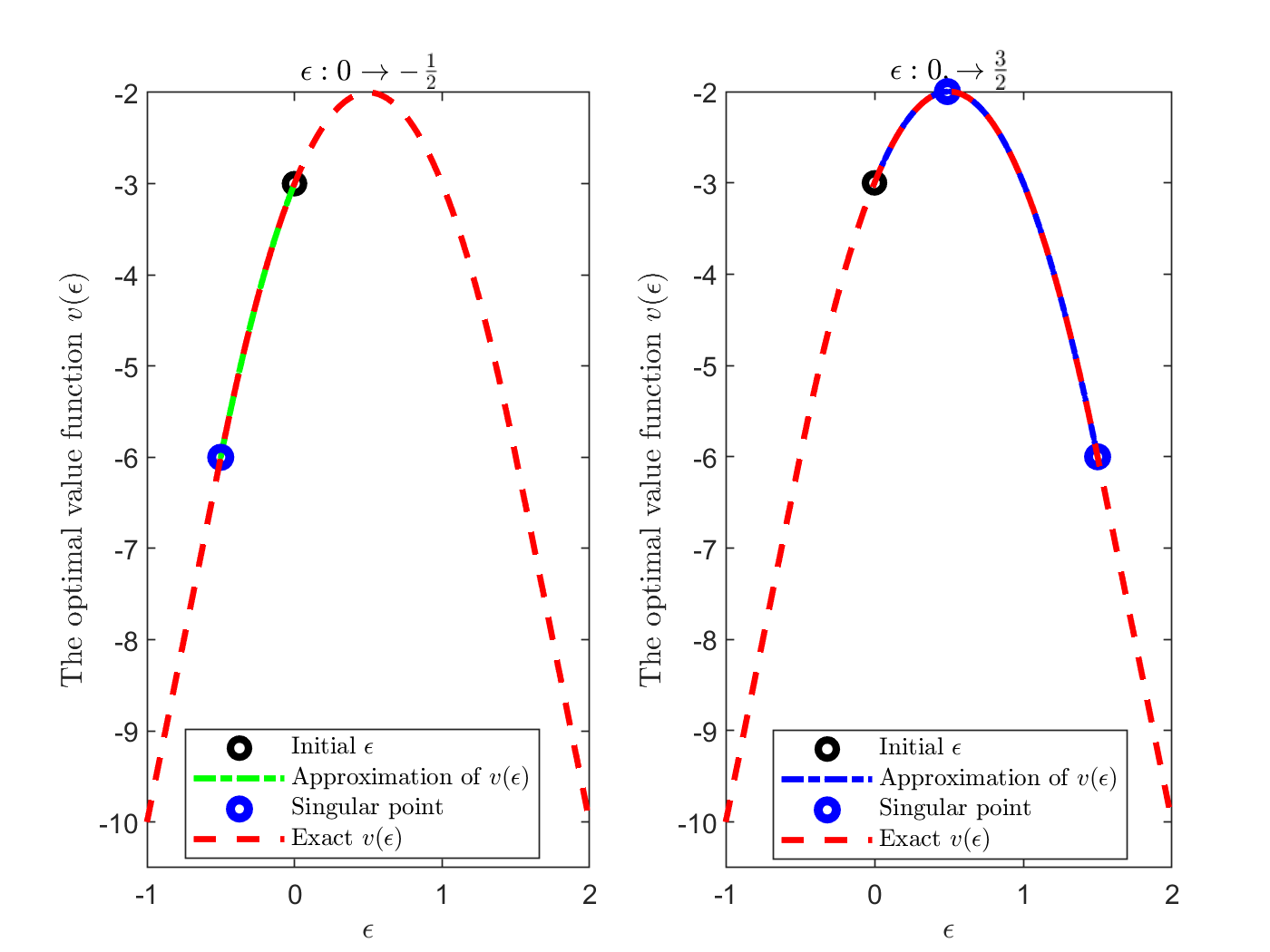}
\end{minipage}\hfill
\begin{minipage}[c]{0.4\textwidth}
\caption{The exact and numerical approximation of 
the optimal value function for problem~\eqref{3elliptope_cut_eq} on $[-1,2]$.}
\label{Fig:elliptope_optimal_val}
\end{minipage}
\end{figure}

Here, we apply Algorithm~\ref{outermost_alg} to identify the singular points and the transition points of the parametric SDO problem~\eqref{3elliptope_cut_eq} in a bounded domain $[-1,2]$. We initialize Algorithm~\ref{outermost_alg} with the initial point $\epsilon = 0$ and the initial increment change $\Delta\epsilon=0.005$. While tracking forwards, Algorithm~\ref{singular_identification} computes the numerical approximation of the unique optimal solution until it locates the singular points $\epsilon=\frac12$ and $\epsilon=\frac32$. Then restarting the solution tracking at $\frac32+\Delta\epsilon$, Algorithm~\ref{invariancy_identification} identifies the invariancy interval $(\frac32,2)$ and the transition point $\epsilon=\frac32$. In an analogous fashion, while tracking backwards, Algorithm~\ref{singular_identification} and Algorithm~\ref{invariancy_identification} identify the singular point $\epsilon = -\frac12$ and the invariancy interval $(-1,-\frac12)$, respectively. Figure~\ref{Fig:elliptope_optimal_val} illustrates the exact and numerical approximation of the optimal value function.

\vspace{5px}
\noindent
Applying Algorithm~\ref{transition_identification} to the singular point $\epsilon=\frac{1}{2}$, we can observe that $V^a(\tfrac12)$ is not isolated, and it belongs to a 1-dimensional irreducible component of $\mathbf{V}\big(F(V,\tfrac12)\big)$. We then invoke the polynomial solver {\tt Bertini} to compute a generic solution
\begin{align*}
X^*(\tfrac12)=\begin{pmatrix}
1 &-0.0449 &-0.0449& 0\\
-0.0449&1 &1& 0 \\
-0.0449& 1 & 1 & 0\\
0 &0&0& 0.0898
\end{pmatrix}, \quad
 y^*(\tfrac12)=\begin{pmatrix} \ \ 0\\ -1\\ -1\\ \ \ 0\\ \ \ 0\\ \ \ 0\\ \ \ 0 \end{pmatrix}, \quad
S^*(\tfrac12)=\begin{pmatrix}
0& \ \ 0&\ \ 0&0\\
0& \ \ 1&-1&0\\
0&-1&\ \ 1&0\\
0& \ \ 0&\ \ 0&0
\end{pmatrix}, 
\end{align*}
in which $\rank\big(X^*(\tfrac12)\big)=3$ and $\rank\big(S^*(\tfrac12)\big)=1$. Given the rank of $X^*(\epsilon)$ and $S^*(\epsilon)$ on $(-\frac12,\frac12) \cup (\frac12,\frac32)$, all this implies that the singular point $\epsilon = \frac12$ belongs to the nonlinearity interval $(-\frac12,\frac32)$. Consequently, the domain $(-1,2)$ is partitioned as 
\begin{align*}
\mathcal{U}_{\mathrm{inv}}=(-1,-\tfrac12) \cup (\tfrac32, 2), \quad \mathcal{U}_{\mathrm{non}}=(-\tfrac12,\tfrac32), \quad \mathcal{U}_{\mathrm{tran}}=\{-\tfrac12,\tfrac32\}.
\end{align*}
	
\section{Concluding remarks and future research}\label{conclusion}
This paper utilized an optimal partition approach for the parametric analysis of SDO problems, where the objective function is perturbed along a fixed direction. In terms of continuity, we provided sufficient conditions for the existence of nonlinearity intervals. Furthermore, we invoked the semi-algebraicity of the optimal set to prove the finiteness of the set of transition points. We showed that the optimal set mapping might fail to be continuous on a nonlinearity interval, and the sequence of maximally complementary optimal solutions may converge to the boundary of the optimal set at an $\epsilon$ in a nonlinearity interval. Finally, under the local nonsingularity condition of Theorem~\ref{thm:V_ODE}, we developed Algorithms~\ref{singular_identification} and~\ref{transition_identification} to compute nonlinearity intervals and identify transition points in $\interior(\mathcal{E})$. If we further assume the \textcolor{black}{generic} global nonsingularity condition of Proposition~\ref{extrema2singular}, Algorithm~\ref{outermost_alg} efficiently partitions $\interior(\mathcal{E})$ into finite union of invariancy intervals, nonlinearity intervals, and transition points. The computational approach was demonstrated on several examples.

\vspace{5px}
\noindent
It is worth mentioning that our optimal partition approach is particularly useful in the context of reoptimization of SDO problems, e.g., matrix completion problems, when the maximal rank of optimal solutions is concerned. Given the lack of efficient warm-start procedures for IPMs, our approach avoids the need for reapplying IPMs after a small perturbation to the objective function, if the given $\epsilon$ belongs to a nonlinearity interval. We should note, however, that quadratic convergence of IPMs is impaired by the failure of strict complementarity or non-degeneracy conditions~\cite{AHO98}, which is always the case at a transition point. Therefore, it would be also interesting to see how computational complexity of IPMs varies on the closure of nonlinearity intervals, e.g., when $\epsilon$ is perturbed from/to a transition point to/from a point in a nonlinearity interval. This is in fact the continuation of the work in~\cite[Section~4]{MT20}, where we provided bounds on the distance between central solutions and approximations of the optimal partitions of the original and perturbed SDO problems.

\vspace{5px}
\noindent
We conjecture that condition~\eqref{continuity_selection} could fail 
at a boundary point of a nonlinearity interval. 
It is worth providing a counterexample or sufficient conditions which guarantee the validity of~\eqref{continuity_selection} at a boundary point of a nonlinearity interval. Furthermore, we still do not know whether the subspaces $\big(\mathcal{B}(\epsilon),\mathcal{T}(\epsilon),\mathcal{N}(\epsilon)\big)$ vary continuously on a nonlinearity interval. These topics are subjects of future research.

\section*{Acknowledgments.}
We are indebted to the anonymous
referees whose insightful comments helped us improve the
presentation of this paper. The first and third authors were supported in part by
Office of Naval Research (ONR) grant N00014-16-1-2722
and National Science Foundation (NSF) grant CCF-1812746.
The second and fourth authors were 
supported by the Air Force Office of Scientific Research (AFOSR) grant FA9550-15-1-0222.

\begin{APPENDIX}{Proofs of Theorems}
\proof{Proof of Theorem~\ref{transition_finite}.} 
Recall that given $\bar{\epsilon} \in \interior(\mathcal{E})$ and a maximally complementary optimal solution $\big(X^*(\bar{\epsilon}),y^*(\bar{\epsilon}),S^*(\bar{\epsilon})\big)$, the ranks of $X^*(\bar{\epsilon})$ and $S^*(\bar{\epsilon})$ are maximal on $\mathcal{P}^*(\bar{\epsilon}) \times \mathcal{D}^*(\bar{\epsilon})$. Hence, the set of all $\epsilon$ with an optimal partition associated with a fixed rank $(\theta,\sigma)$ can be defined as 
\begin{align*}
\mathcal{S}_{(\theta,\sigma)}\!:=\!\Big\{\epsilon \in \mathbb{R} : \exists \ (X,y,S) \in \ri\!\big(\mathcal{P}^*(\epsilon) \times \mathcal{D}^*(\epsilon)\big), \ \rank(X)=\theta, \ \rank(S)=\sigma, \ \epsilon \in \interior(\mathcal{E}) \Big\},
\end{align*} 
which in turn implies 
\begin{align}\label{full_coverage}
\interior(\mathcal{E})=\bigcup_{\substack{\theta,\sigma\in\{0,\ldots, n\}\\ \theta+\sigma \le n}} \mathcal{S}_{(\theta,\sigma)},
\end{align}
\textcolor{black}{where $(\theta,\sigma)$ is a pair of integers}. In what follows, we prove that $\bar{\epsilon}$ is a transition point if and only if $\bar{\epsilon} \in \bd(\mathcal{S}_{(\theta,\sigma)}) \cap \interior(\mathcal{E})$ for some nonnegative integer $(\theta,\sigma)$ with $\theta+\sigma \le n$, and that $\mathcal{S}_{(\theta,\sigma)}$ is a semi-algebraic subset of $\mathbb{R}$. Then the finiteness follows from the fact that $\mathcal{S}_{(\theta,\sigma)}$ has only a finite number of boundary points~\cite[Theorem~5.22]{BPR06}.

\paragraph{Equivalency of boundary points and transition points} By Definition~\ref{def:transition_point}, it is clear that if $\hat{\epsilon} \in \interior(\mathcal{E})$ is a boundary point of $\mathcal{S}_{(\theta,\sigma)}$, then $\hat{\epsilon}$ must be a transition point. More specifically, by the definition of a boundary point,
\begin{itemize}
\item if $\hat{\epsilon} \not \in \mathcal{S}_{(\theta,\sigma)}$, then every neighborhood of $\hat{\epsilon}$ contains an $\epsilon' \in \mathcal{S}_{(\theta,\sigma)}$, which implies that either $\rank\!\big(X^*(\epsilon')\big) \neq \rank\!\big(X^*(\hat{\epsilon})\big) $, $\rank\!\big(S^*(\epsilon')\big) \neq \rank\big(S^*(\hat{\epsilon})\big)$, or both holds;
\item if $\hat{\epsilon} \in \mathcal{S}_{(\theta,\sigma)}$, then every neighborhood of $\hat{\epsilon}$ contains an $\epsilon'' \in \interior(\mathcal{E}) \setminus \mathcal{S}_{(\theta,\sigma)}$, which implies that either $\rank\!\big(X^*(\epsilon'')\big) \neq \rank\!\big(X^*(\hat{\epsilon})\big) $, $\rank\!\big(S^*(\epsilon'')\big) \neq \rank\!\big(S^*(\hat{\epsilon})\big)$, or both holds.
\end{itemize}
From both cases, it is immediate that $\hat{\epsilon}$ is a transition point. Conversely, by~\eqref{full_coverage}, a transition point $\bar{\epsilon}$ belongs to $\mathcal{S}_{(\theta,\sigma)}$ for some nonnegative integer $(\theta,\sigma)$ with $\theta+\sigma \le n$. If $\bar{\epsilon} \in \interior(\mathcal{S}_{(\theta,\sigma)})$, then the ranks of $X^*(\epsilon)$ and $S^*(\epsilon)$ would be constant on a neighborhood of $\bar{\epsilon}$, which is a contradiction. Therefore, we must have $\bar{\epsilon} \in \bd(\mathcal{S}_{(\theta,\sigma)})$, see e.g.,\textcolor{black}{~\cite[Page~102]{M2000}}, which completes the first part of the proof.

\paragraph{Semi-algebraicity of $\mathcal{S}_{(\theta,\sigma)}$}
We proceed with the proof of semi-algebraicity in three steps. For the ease of exposition and by using the isometry~\eqref{isomorphism}, we sometimes identify the optimal solutions by column vectors $V=(x; y; s)$, \textcolor{black}{where $x$ and $s$ are obtained from the upper triangular entries of $X$ and $S$, respectively}. 
\begin{enumerate}
\item Given a fixed $\epsilon$, $\mathcal{P}^*(\epsilon) \times \mathcal{D}^*(\epsilon)$ is the set of all vectors $V$ satisfying~\eqref{KKT_conditions} and~\eqref{eq:KKT_INEQ}, where~\eqref{eq:KKT_INEQ} is equivalent to $2(2^n-1)$ polynomial inequalities, enforcing all principal minors of $X$ and $S$ to be nonnegative. Therefore, $\mathcal{P}^*(\epsilon) \times \mathcal{D}^*(\epsilon)$ is a semi-algebraic subset of $\mathbb{R}^{m+2t(n)}$, i.e., $\mathcal{P}^*(\epsilon) \times \mathcal{D}^*(\epsilon)$ is defined by a Boolean combination of polynomial equalities and inequalities~\cite[Page~57]{BPR06}.
\item Since $\mathcal{P}^*(\epsilon) \times \mathcal{D}^*(\epsilon)$ is convex, see e.g.,~\cite[Theorem~6.4]{Rock70}, the relative interior of $\mathcal{P}^*(\epsilon) \times \mathcal{D}^*(\epsilon)$ is the set of all $V$ satisfying
\begin{align*}
\forall \ \bar{V} \in \mathcal{P}^*(\epsilon) \times \mathcal{D}^*(\epsilon), \quad  \exists \ \gamma > 0 \quad \st \quad  V + \gamma (V-\bar{V}) \in \mathcal{P}^*(\epsilon) \times \mathcal{D}^*(\epsilon), 
\end{align*}
which, by semi-algebraicity of $\mathcal{P}^*(\epsilon) \times \mathcal{D}^*(\epsilon)$, can be expressed by a quantified formula $\Psi$ (a formula with quantifiers from the set $\{\forall,\exists\}$) in the language of ordered fields, see e.g.,~\cite[Proposition~3.1]{BPR06}. A \textit{formula}~\cite[Page~13]{BPR06} is the Boolean combination of polynomial equalities and inequalities with real coefficients. Since the $\mathbb{R}$-realization of $\Psi$, i.e., the set of all real solutions satisfying $\Psi$, is a semi-algebraic subset of $\mathbb{R}^{m+2t(n)}$~\cite[Theorem~2.77]{BPR06}, we just showed that $\ri\!\big(\mathcal{P}^*(\epsilon) \times \mathcal{D}^*(\epsilon)\big)$ is also a semi-algebraic subset of $\mathbb{R}^{m+2t(n)}$. \label{relative_interior_semi-algebraic}
\item The set $\{x \in \mathbb{R}^{t(n)} : \rank(x) = \theta\}$ is equal to
\begin{align*}
\{x \in \mathbb{R}^{t(n)} : \rank(x) = \theta\}= \{x \in \mathbb{R}^{t(n)} : \rank(x) \le \theta\} \cap \big(\mathbb{R}^{t(n)} \setminus \{x \in \mathbb{R}^{t(n)} : \rank(x) \le \theta-1\}\big),
\end{align*}
where $\{x \in \mathbb{R}^{t(n)} : \rank(x) \le \theta\} = \{x \in \mathbb{R}^{t(n)} : \text{all minors of $x$ of size $\theta+1$ are zero}\}$, see e.g.,~\cite[Page~12]{HJ12}, is an algebraic set, as minors of $x$ are polynomials in terms of the entries of $x$. This also implies that $\mathbb{R}^{t(n)} \setminus \{x \in \mathbb{R}^{t(n)} : \rank(x) \le \theta-1\}$ is a semi-algebraic subset of $\mathbb{R}^{t(n)}$~\cite[Page~57]{BPR06}. \label{rank_semi-algebraic}
\end{enumerate}

\vspace{10px}
Using the arguments in~\eqref{relative_interior_semi-algebraic} and~\eqref{rank_semi-algebraic}, and given a fixed $(\theta,\sigma)$, the set 
\begin{align}\label{semi-algebraic_fixed_rank}
\!\big\{(V,\epsilon) \in \mathbb{R}^{m+2t(n)+1} :  V \in \ri\!\big(\mathcal{P}^*(\epsilon) \times \mathcal{D}^*(\epsilon)\big), \ \rank(x)=\theta, \ \rank(s)=\sigma,  \ \epsilon \in \interior(\mathcal{E})  \big\}
\end{align} 
is a semi-algebraic subset of $\mathbb{R}^{m+2t(n)+1}$, because it is the $\mathbb{R}$-realization of a quantified formula. As a result, the projection of~\eqref{semi-algebraic_fixed_rank} to $\mathbb{R}$, i.e., $\mathcal{S}_{(\theta,\sigma)}$ is a semi-algebraic subset of $\mathbb{R}$~\cite[Theorem~2.76]{BPR06}, which completes the second part of the proof.
 \halmos
\endproof
\vspace{10px}

\proof{Proof of Theorem~\ref{thm:V_ODE}.}
By Lemma~\ref{nonsingularity_conditions}, $V^*(\epsilon)$ is the unique optimal solution of $(\mathrm{P_{\epsilon}})-(\mathrm{D_{\epsilon}})$ with a nonsingular Jacobian for every $\epsilon \in \mathcal{I}_{reg}$. Thus, by the analytic implicit function theorem~\cite[Theorem~10.2.4]{D60}, $V^*(\epsilon)$ is analytic on $\mathcal{I}_{reg}$. On the other hand, since $V^*(\epsilon)$ satisfies~\eqref{KKT_conditions} point-wise, it is easy to see, by taking the derivatives of the equations in~\eqref{KKT_conditions}, that $V^*(\epsilon)$ is an analytic solution of the ODE system~\eqref{ODE_system}.

\vspace{5px}
\noindent
Now, let us consider a differentiable mapping $\underline{V}(\epsilon)\!:=\big(\underline{X}(\epsilon),\underline{y}(\epsilon),\underline{S}(\epsilon)\big)$ as an arbitrary solution of~\eqref{ODE_system}. Then $\underline{V}(\epsilon)$ solves~\eqref{KKT_conditions} point-wise, $\underline{V}(\bar{\epsilon})=V^*(\bar{\epsilon})$, and $J(\underline{V}(\epsilon),\epsilon)$ is nonsingular on $\mathcal{I}_{reg}$, because the right hand side of~\eqref{ODE_system} must be bounded on $\mathcal{I}_{reg}$. By invoking the nonsingularity of $J\!\big(V^*(\bar{\epsilon}),\bar{\epsilon}\big)$, and using the analytic implicit function theorem, we can immediately see that $\underline{V}(\epsilon)=V^*(\epsilon)$ on a neighborhood of $\bar{\epsilon}$. However, if we further take into account the nonsingularity of $J\!\big(\underline{V}(\epsilon),\epsilon\big)$ on $\mathcal{I}_{reg}$ and apply the analytic implicit function theorem again, then $\underline{V}(\epsilon)$ must be analytic on $\mathcal{I}_{reg}$ as well. Therefore, as a result of~\cite[Corollary~1.2.6]{KP02}, $\underline{V}(\epsilon)=V^*(\epsilon)$ holds globally on $\mathcal{I}_{reg}$. This completes the proof of uniqueness of $V^*(\epsilon)$.  \halmos
\endproof

\end{APPENDIX}


\bibliographystyle{siam}   
\bibliography{mybibfile}

\end{document}